\documentclass[10pt,twoside]{article}
\usepackage{amssymb}
\usepackage{latexsym}
\usepackage{fullpage}
\usepackage{latexsym,rotating,amsmath,epsfig}
\usepackage{here}
\usepackage{amsmath}
\usepackage{amscd}
\usepackage{oldgerm}
\usepackage{fancyheadings}
\def\In{\mathbb{I}}
\def\N{\mathbb{N}}

\def\R{\mathbb{R}}

\def\Bb{\overline{B}}
\def\1{\mathbf{1}}
\def\:{\lrcorner}
\def\#{\sharp}
\def\k{\kappa}

\def\l{\lambda}
\def\a{\alpha}
\def\b{\beta}
\def\g{\gamma}
\def\d{\delta}
\def\e{\epsilon}

\def\o{\circ}

\def\s{\sigma}

\def\x{\otimes}
\def\<#1,#2>{\langle#1,\,#2\rangle}

\def\G{\Gamma}

\def\S{\mathbb{S}\,}

\def\qed{\ensuremath{\quad\Box\quad}}
\def\pfill{\par\vskip2mm plus1mm minus1mm\noindent}
\def\inv#1{\raise.1em\hbox to 0pt{$^{-1}$\hss}_{#1}\;}

\def\v{\noindent}

\newcommand{\bean}{\begin{eqnarray*}}
\newcommand{\eean}{\end{eqnarray*}}
\newcommand{\benu}{\begin{enumerate}}
\newcommand{\eenu}{\end{enumerate}}
\newcommand{\eea}{\end{eqnarray}}
\newcommand{\bea}{\begin{eqnarray}}

\newtheorem{Theorem}{Theorem}[section]
\newtheorem{Lemma}[Theorem]{Lemma}
\newtheorem{Corollary}[Theorem]{Corollary}
\newtheorem{Definition}[Theorem]{Definition}

\newtheorem{Proposition}[Theorem]{Proposition}

\title{A metric approach to Fr\'echet geometry}
\begin{document}

\author{Olaf M\"uller\footnote{Instituto de Matem\'aticas, Universidad Nacional Aut\'onoma de M\'exico (UNAM) Campus Morelia, C. P. 58190, Morelia, Michoac\'an, Mexico. email: olaf@matmor.unam.mx}}

\maketitle

\begin{abstract}
\v The aim of this article is to present the category of bounded Fr\'echet manifolds in respect to which we will review the geometry of Fr\'echet manifolds with a stronger accent on its metric aspect. An inverse function theorem in the sense of Nash and Moser in this category is proved, and some applications to Riemannian geometry are given. 
\end{abstract}

\section{Introduction}

\v The celebrated Nash-Moser inverse function theorem is one of the most powerful tools linking infinite-dimensional geometry with global analysis when dealing with smooth objects lacking a Banach structure. However, the usual description of the tame category uses the notion of a series of seminorms which is intuitively not very well-accessible, also it delivers not much geometric insight in these infinite-dimensional spaces (in most cases spaces of smooth sections). So an idea to remedy this gap could be to fix natural metrics on these spaces of sections which can encode the crucial properties of tameness. This will be done in the following article. Instead of tame maps (a notion that will refer to a metric instead of to a series of seminorms) we will often as well consider maps which are bounded w.r.t. the metrics used. For maps of finite differentiability this procedure would not be well-applicable as fundamental operators like covariant derivatives are unbounded, meanwhile in spaces of {\em smooth} maps with the usual metrics they {\em are} bounded. One reason to consider metrics instead of series of seminorms is certainly that they are geometrucally more demonstrative. A second reason to be interested in concrete metrics on these spaces comes from applications in the geometric quantization of field theories. A third advantage consists in the existence of many fix point theorems for metric spaces. Thus the basic question of this article is: {\em How metric can an exhaustive description of the properties of Fr\'echet spaces be?}

\v The article, apart of presenting new results, is thought to serve as an introduction to Fr\'echet spaces from the metric point of view; in order to keep it as self-contained as possible, several well-known results (indicated as such) were included as well. It is structured as follows: In section 2 we review some basic results for connections in fiber bundles. The treatment of vector bundle connections as Levi-Civita connections on the total space does not seem to be mentionned in the literature so far. Likewise, some technical definitions and results about bounded geometry are established. In section 3, the category of bounded Fr\'echet manifolds is presented and applied to spaces of sections of fiber bundles. Finally, in section 4 we prove an inverse function theorem fitting to this category, right along the lines of the proof of the inverse function theorem for maps between Banach space. The last sections treats the question whether the notion of length structure is applicable in the present context.  

\v I would like to express my gratitude towards Helge Gl\"ockner who drew my attention to a crucial mistake in the first versions of this article and made me familiar with a special class of metrics.

\section{Preliminaries: Connections in fiber bundles and bounded geometry} 

\subsection{General fiber bundles}

\v As a preparation, we review some facts from the theory of connections.

\begin{Definition}
Let a fibre bundle $\pi : N\rightarrow M$ be given, $m= dim (M)$. A {\bf connection on $\pi$} is a smooth $m$-dimensional horizontal distribution in $\tau_N$, i.e. an $m$-dimensional subbundle $D$ of $\tau_N$ with $ker (d \pi \vert_{D} ) = \{ 0 \}$ at each point.
\end{Definition}

\begin{Definition}
Let $\pi: N \rightarrow M$ be a fibre bundle with a connection $D$. Let a curve $c: [0,1] \rightarrow M$ be given. Let $p \in \pi ^{-1} (c(0))$ be given. Then a {\bf horizontal lift of $c$ starting at $p$} is the curve $L^D_{c,p} : [0,1] \rightarrow N$ which is given by $ L^D_{c,p} (0) := p$ and $\dot{L}^D_{c,p} (t) := (d \pi \vert_{D}) ^{-1} (\dot{c} (t))$. By definition, $L^D_{c,p}$ projects to $c$ under $\pi$.
\end{Definition}

\v Horizontal lifts are unique, but they do not exist always. If a connection gives rise to a horizontal lift for every curve in the base manifold, the connection is said to be {\bf Ehresmann-complete} or an {\bf Ehresmann-connection}.

\v It is easy to see by considering inverse curves that in the case of an Ehresmann connection each curve $c: [0,1] \rightarrow M$ in the base manifold induces a diffeomorphism $D_c$ from $\pi^{-1} (c(0))$ to $\pi^{-1} (c(1))$ which we call {\bf parallel transport along $c$}. In the same manner, a curve of diffeomorphisms $H_{\l}$ of $M$ starting at the identity induces a curve of fiber bundle isomorphisms $D^H_{\l}$. If we have a metric fiber bundle with a connection for which the parallel transports along all curves are isometries of the fibers we call the connection {\bf metric}.

\begin{Definition}[Pull-back of connections by immersions]
Let $\pi: N \rightarrow M$ be a fiber bundle with a connection $D$. Let $f:S \rightarrow M$ be an immersion of a manifold $S$ into $M$. Then we can equip the pull-back bundle $f^* \pi$ with the {\bf pull-back connection} which is denoted by $f ^* D := T(f^* N) \cap D$ where $f^* N$ is understood as a submanifold of $N$.  
\end{Definition}

\v It is easy to check that this is really a connection: First note that because of the finite-dimensional inverse function theorem around each point $p \in S$ there is a neighborhood $U$ whose image is a submanifold of $M$. Then, by considering a locally trivial submanifold chart neighborhood of $f(U) \subset M$ we show that $f^* N = \pi ^{-1} (f(U))$ is a submanifold. Then, for $\k$ being the submanifold chart and for $t$ being the trivialization, by considering the chart $(\k \times \1) \o t $ for $\pi ^{-1} (U)$ we see that $T(f^* N) = \{ v \in TN \vert d \pi (v) \in T(f(U)) \} $. As $rang (D) = m$ and $corang (T(f^*N))= m-s$, we get $rang(f^* D) \geq s $. On the other hand, $d \pi (v) \in T(f(U))$ for $v \in T(f^* N)$, and $d \pi (v) \neq 0$ for $v \in (D \cap \pi^{-1} (m)) \setminus \{  0 \} $. Therefore $ker(d (f^* \pi ) \vert _{f^* D} = \{ 0 \} $, and the claim follows.    

\begin{Proposition}
\label{pullback}
Between the parallel transports of a connection and its pull-back connection under an immersion $f$ we have the simple relation $f^*D_c = D_{f\o c}$.
\end{Proposition}

\begin{Definition}
Let $\pi: N \rightarrow M$ be a fiber bundle with a connection $D$. A vector field $X$ on $N$ is called {\bf projectable} iff $\pi_* (X(p)) = \pi_* (X(q))$ for all $p,q \in \pi ^{-1} (m)$. It is called {\bf basic} if it is projectable and horizontal. 
\end{Definition}

\begin{Proposition}
\label{Lie-fiber}
Let $\pi$ be a fiber bundle with a connection. If $X$ is a projectable vector field and $Y$ is a vertical vector field, then $[X,Y]$ is again vertical.
\end{Proposition}

\v {\bf Proof.} This follows from $[X,Y]= \mathcal{L}_X Y = \frac{d}{dt} \vert_{t=0} (Fl_{-t}^X \o Y \o Fl_{t}^X ) $ and from the fact that the local flow of $X$ (which is a lift of the local flow of $\pi_* X$ and hence is defined on the whole fiber) translates fibers into fibers. \hfill \qed

\bigskip

\v Let a fiber bundle $\pi: E \rightarrow M$ be given with a pseudo-Riemannian metric $h$ on the fibres, a pseudo-Riemannian metric $g$ on the base manifold $M$, and a connection $D$. Then we can use $D$ to equip the total space with the {\bf total metric} $G=G_{g,h,D}$ which makes the projection a pseudo-Riemannian submersion. Then we can consider the Levi-Civita connection $\nabla^E$ on the total space.

\v Now in the case that the metrics are Riemannian, the total Riemannian metric $G_{g,h,D}$ provides us with convenient criteria to determine whether a given connection is Ehresmann:

\begin{Theorem}
Let $\pi: E \rightarrow M$ be a bundle with a connection $D$. If, for some appropriately chosen Riemannian metrics $g,h$ on $M$ resp. on the fibers, $(E, G_{g,h,D})$ is complete, then $D$ is Ehresmann. Likewise, if we can find complete metrics $g$ and $h$ such that the distortion of $D_c$ as a map between the fibers is locally bounded with respect to the fiber metrics $h$, then $D$ is Ehresmann.  \end{Theorem}

\v {\bf Proof} in \cite{be}, 9.42 and 9.46 \hfill \qed

\subsection{Vector bundle connections}

\v For vector bundle connections we require additionally to the previous that they respect the linear structure: 

\begin{Definition}(cf. \cite{be}, 9.51  \footnote{\v Note, however, that contrary to the claim of this otherwise excellent reference it is not sufficient that the stucture group be {\em conjugate} to $Gl(n)$; it has to be $Gl(n)$ itself which is a distinguished subgroup of the diffeomorphisms in the fibre as soon as we have a vector bundle.})
Let $\pi:E \rightarrow M$ be a vector bundle. A {\bf vector bundle connection on $\pi$} is a connection on $\pi$ as a fiber bundle whose parallel transports are {\em linear} maps of the corresponding fibers.   
\end{Definition}

\v {\bf Remark.} In particular the definition above says that the parallel transport is defined overall, i.e. that the connection is Ehresmann.

\bigskip

\v A vector bundle connection gives rise to a covariant derivative in the following way: First, vector addition provides us with a map $\iota: T_p E \rightarrow \pi ^{-1} (m)$ for all $p \in \pi ^{-1} (m)$. Then for a vector $X \in T_p M$ and a section $\g$ of $\pi: E \rightarrow M$ its covariant derivative $\nabla_X \g $ is the section of $\pi $ defined by  

$$ \nabla_X \g (p) := \iota (\frac{d}{ds} (D^{-1}_{c[s]} (\gamma \o c(s)) )\vert_{s=0} ) =\iota ( \lim_{s \rightarrow 0} \frac{D^{-1}_{c[s]} ( \g (c(s))) - \g (c (0))}{s}) = \iota ( \lim_{s \rightarrow 0} \frac{ \g (c(s)) -  D_{c[s]} ( \g (c (0)))}{s})$$

\v in which $c$ is a curve with $c(0) = (p,t) $, $\dot{c} (0) = X$ and $c[s]:= c \vert_{[0,s]}$. One immediately checks that this independent of the choice of $c$.

\begin{Proposition}
Let $D$ be a metric connection on the Pseudo-Riemannian vector bundle $(\pi:E \rightarrow M, g)$, let $f: S \rightarrow M$ be an embedding, then $f^* D$ is metric on $(f^* \pi, f^* g)$. 
\end{Proposition}

\v {\bf Proof} by Proposition \ref{pullback}  \hfill \qed

\bigskip

\v We have a correspondance $\tilde{}$ (inverse to $\iota$ as above) between a point $v$ in the total space and a vertical translational-invariant vector field $\tilde{v}$ on the fiber. Translational-invariant vector fields are always vertical. Translational-invariant vector fields form a maximal center in the algebra of vector fields in the total space.

\begin{Proposition}
\label{Lie}
Let $\pi: N \rightarrow M$ be a vector bundle with vector bundle connection and arbitrary (not necessarily translational-invariant) metric on the fibers, let $s$ be a section of $\pi$ which is parallel along a vector field $X$ on $M$, let $\hat{X}$ the associated basic vector field to $X$, then $[\hat {X}, \tilde{s}]$ is translational-invariant. 
\end{Proposition}

\v {\bf Proof.} This is because the finitesimal commutator of the flow of $\hat{X}$ (which is parallel transport along the flow of $X$ and hence linear) and the flow of $\tilde{s}$ (a translation in every fiber) is an affine map. \hfill \qed

\bigskip

\v In the following we examine the relation between a vector bundle connection and the Levi-Civita connection on the total space to the metric $G_{g,h,D}$ as in the previous section.

\begin{Proposition}
\label{vertigo}
Let $\pi$ be a fiber bundle with a Riemannian metric on the base space and the fibers and with a connection $D$ as above, let $\nabla^E$ denote the Levi-Civita connection as above. Let $X$ be a horizontal vector field and $Y$ a vertical vector field. Then $\nabla^E_X Y $ is vertical.
\end{Proposition}

\v {\bf Proof.} Consider one point $p$ in the total space. Because of tensoriality, to calculate $\nabla_X Y$ at $p$, we can assume that $X$ is projectable. Let $V$ be another projectable vector field and show $\langle \nabla^E _X Y, V \rangle =0 $ by using Lemma \ref{Lie-fiber}. \hfill \qed

\begin{Proposition}
Let a vector bundle $\pi$ be given with a vector bundle connection and a translational-invariant metric. Let $X$ be basic and $V$ be translational-invariant. Then $\nabla_X V$ is translational-invariant.
\end{Proposition}

\v {\bf Proof.} First recall that, as $V$ is vertical and because of Proposition \ref{vertigo}, $\nabla_X V$ is vertical. We show that it is invariant under {\em horizontal} translations, that means, under addition of vector fields of the form $\tilde{Y}$ where $Y$ is a horizontal section of $\pi$. Now let $\tilde{Y}$ be as above and $Z$ be another translational-invariant vector field, then

\begin{align*}
\tilde{Y} (\langle \nabla_X V, Z \rangle) &= \frac{1}{2} \tilde{Y} ( X \langle V, Z \rangle - \langle X, [V,Z] \rangle + \langle Z, [X,V] \rangle + \langle V, [Z,X] \rangle )  \\
&= \frac{1}{2} (X \tilde{Y} \langle V, Z \rangle - [X, \tilde{Y}] -0 + \langle Z, [\tilde{Y}, [X,V]] \rangle + \langle V, [\tilde{Y}, [X,Z]] \rangle)  \\ 
&=  \frac{1}{2}( \langle Z, [\tilde{Y}, [X,V]] \rangle + \langle V, [\tilde{Y}, [X,Z]] \rangle)  \\
&= \frac{1}{2} (\langle - Z, [X, - [V, \tilde{Y}]] - [V, [\tilde{Y}, X]] \rangle + \langle V , -[Z, [X, \tilde{Y}]] - [X, [\tilde{Y}, Z]] \rangle )  \\ 
&=0  
\end{align*}

\v where the first step is the Koszul formula, the second step is the definition of the Lie derivative resp. the translational invariance of the fiber metric, the third step uses Proposition \ref{Lie-fiber}, the fourth step is the Jacobi identity, and the last step uses Proposition \ref{Lie-fiber} again and the fact that translations commute. \hfill \qed

\begin{Proposition}
\label{parallel}
Let a vector bundle $\pi$ be given with a metric vector bundle connection and a translational-invariant metric.
Let $\hat{X}$ be the basic vector field on $E$ to a vector field $X$ on $M$. Let $Y$ be an $X$-parallel section of $\pi$, then $\nabla_{\hat{X}} \tilde{Y} =0$. \end{Proposition}

\v {\bf Proof.} Choose another translational-invariant vector field $\hat{Z}$, w.r.o.g. $\nabla_X Z =0$ at the point $p$ in question, then the Koszul formula gives at $p$

\begin{align*}
\langle \nabla_{\hat{X}} \tilde{Y}, \hat{Z} \rangle &= \hat{X} \langle \tilde{Y}, \tilde{Z} \rangle - \langle \hat{X} , [\tilde{Y}, \tilde{Z}] \rangle + \langle \tilde{Z}, [ \hat{X}, \tilde{Y}] \rangle + \langle \tilde{Y}, [ \hat{X}, \tilde{Z}] \rangle \\
&= X \langle Y, Z \rangle + \langle \tilde{Z},   \lim_{t \rightarrow 0} ( P_t^X \o \tilde{tY}- \tilde{tY} \o P_t^X)  \rangle   + \langle \tilde{Y},   \lim_{t \rightarrow 0} (P_t^X \o \tilde{tZ}- \tilde{tZ} \o P_t^X ) \rangle\\
&= \langle \tilde{Z},   \lim_{t \rightarrow 0} \widetilde{P_t^X (Y)}   \rangle   + \langle \tilde{Y},   \lim_{t \rightarrow 0} \widetilde{P_t^X ( Z)}  \rangle =0
\end{align*}

\v where the second equation is due to the fact that $\langle \tilde{Y}, \tilde{Z} \rangle$ constant on the fiber and the last equation is due to the fact that parallel transport of a metric connection is an orthogonal map and its derivative is a skew-symmetric map.    \hfill \qed

\bigskip

\v Now we can relate the connections by the following

\begin{Theorem}
\label{huhu}
Let a vector bundle $\pi$ be given with a metric vector bundle connection and a translational-invariant metric. Then we have 
$\widetilde{\nabla^{\pi}_X v } = \nabla^E_{\hat{X}} \tilde{v}  $, where $\hat{X}$ is the basic vector field corresponding to $X$.
\end{Theorem}

\v {\bf Proof.} We have

\begin{align*}
\widetilde{\nabla^{\pi}_X v} &= \lim_{t \rightarrow 0} \widetilde{\frac{P_X^t v - v \o Fl_X^t}{t}}\\
&= \lim_{t \rightarrow 0} \frac{\widetilde{P_X^t v}- \widetilde{v \o Fl_X^t}}{t}\\
&= \lim_{t \rightarrow 0} \frac{P_{\hat{X}}^t \tilde{v} - \tilde{v} \o Fl_{\hat{X}}^t }{t}\\
&= \nabla_{\hat{X}} \tilde{v} 
\end{align*}

\v where the first equation is a consequence continuity of the map $\tilde{}$, the second one of its linearity, the third one of Proposition \ref{parallel}. \hfill \qed

\begin{Theorem}
\label{haha}
Let a vector bundle $\pi$ be given with a vector bundle connection and a translational-invariant metric. Then, if $\tilde{w}$ is translational-invariant and $v$ is vertical, $\nabla_v \tilde{w}$ is horizontal. If additionally the vector bundle connection is metric, then $\nabla_v \tilde{w}=0$.  
\end{Theorem}
 
\v {\bf Proof.} Because of tensoriality we can assume that $v$ is translational-invariant as well. Then the Koszul formula gives immediately $\langle \nabla_v \tilde{w}, z  \rangle =0$ for any vertical vector $z$ which, again by tensoriality, we can assume to be the value of a translational-invariant vector field as well. On the other hand, the application of the same Koszul formula to a horizontal vector $z$ gives

$$ g(\nabla_v \tilde{w} , z ) = \frac{1}{2} (- \mathcal{L} _z (g( v, \tilde{w} ))  + g( v , \mathcal{L} _z \tilde{w} ) + g( \tilde{w}, \mathcal{L}_z v ) = (\mathcal{L}_z g) (v,w)    $$

\v which vanishes if the connection is metric, i.e. if the horizontal flow preserves the metric on the fibers. \hfill \qed

\v Furthermore, there is also a correspondence for basic vector fields:

\begin{Theorem}
\label{hoho}
Let a vector bundle $\pi: E \rightarrow M$ be given with a metric vector bundle connection and a translational-invariant metric. Then for any vector fields $X,Y$ on $M$, for their horizontal lifts $\hat{X}, \hat{Y}$ we have

$$\nabla_{\hat{X}} \hat{Y} = \widehat{\nabla_X Y } .$$

\end{Theorem}

\v {\bf Proof.} This is because $\nabla_{\hat{V}}$ is horizontal: 

$$\langle \nabla_{\hat{V}} \hat{W}, \tilde{X} \rangle = \hat{V} \langle \hat{W}, \tilde{X} \rangle - \langle \hat{W} , \nabla_{\hat{V}} \tilde{X} \rangle = 0$$

\v and because of the corresponding formula controlling the horizontal part in general Riemannian submersions that can be found in \cite{on}, p. 212 (Lemma 45.3) which is proven using the Koszul formula. \hfill \qed 

\bigskip

\v Finally, the last possible combination in the covariant derivative vanishes:

\begin{Theorem}
\label{hihi}
Let a vector bundle $\pi: E \rightarrow M$ be given with a metric vector bundle connection and a translational-invariant metric. Then for any two vector fields $X, V$ on $M$ we have

$$\nabla_{\tilde{X}} \hat{V}=0 .$$

\end{Theorem}

\v {\bf Proof.} First show that $\nabla_{\tilde{X}} \hat{V}$ is horizontal picking a translational-invariant vector field $\tilde{Y}$ and using Theorem \ref{haha}. On the other hand, $\nabla_{\tilde{X}} \hat{V}= \nabla_{\hat{V}} \tilde{X} - [\hat{V}, \tilde{X}] $ must be translational-invariant following the theorems \ref{haha} and \ref{Lie}. This leaves us with the zero vector at every point. \hfill \qed

\bigskip

\v So, for an arbitrary vector bundle $\pi: E \rightarrow M$ with translational-invariant metric and metric vector-bundle connection, and for any curve $c: \In \rightarrow E$ with $k:= \pi \o c$, $\dot{k} \neq 0$, we have

\bea
\label{curve}
\dot{c}(t) = \dot{c}_{vert} (t)+ \dot{c}_{hor} (t), \ \dot{c}_{vert} (t) = \frac{1}{\vert \vert \dot{k} \vert \vert}\widetilde{ \nabla_{k(t)}^{LC} c(t)}
\eea

\v while $\dot{c}_{hor} = \hat{k^. }$. This is because by definition $\nabla_{\hat{k}(t)}^{LC (E)} \tilde{c}(t) =0$ if and only if $\dot{c}$ is horizontal along $\dot{k}$ by Theorem \ref{huhu}.

\subsection{Bounded geometry}
\label{bounded}

\begin{Definition}
Let $F: \pi \rightarrow \pi '$ be a vector bundle homomorphism between two Riemannian vector bundles $\pi, \pi '$ over manifolds $M, N$. The {\bf upper resp. lower distortion $\overline{\d}$ resp. $\underline{\d}$ of $F$} is defined as the quantities

$$ \overline{\d} (F) :=  \sup_{X \in SE} \vert \vert F(X) \vert \vert =  \sup_{X \in E \setminus n} \frac{\vert \vert F( X) \vert \vert}{\vert \vert X \vert \vert},$$

$$ \underline{\d} (F) :=  \sup_{X \in SE} (\vert \vert F(X) \vert \vert^{-1}) \} =  \sup_{X \in E \setminus n} \frac{\vert \vert  X \vert \vert}{\vert \vert F(X) \vert \vert } $$

\v where $S \pi: SE \rightarrow M$ is the sphere bundle of $\pi$ and $n$ its zero section, and the {\bf total distortion $\d$ of $F$} as $\d(F):= max \{ \overline{\d} (F), \underline{\d} (F) \}$.

\end{Definition}

\begin{Definition}
Let $(N,h)$ be a Riemannian manifold, let $V \subset M$ be open. We call $(N,h)$ {\bf of bounded geometry (in $V$)} if for the Levi-Civita connection on $(N,h)$ holds: The injectivity radius $\iota_M$ is bounded from zero on $M$ (resp. on $V$), and $\vert \vert \nabla ^{(k)} R^M (p) \vert \vert_{(\tau^*_M)^3 \x \tau_M} \leq C_k$ for all $k \in \N \cup \{ 0\}$ and for all $p \in M$ (resp. $ p \in V$). 
Let $\pi: E \rightarrow M$ be a fibre bundle $(E,h) \rightarrow (M,g)$ over a Riemannian manifold $(M,g)$ with a fibre bundle connection $D$ and a fiber metric. Let $U \subset E$ be open. We call $\pi$ {\bf of strongly bounded geometry (in $U$)} if the natural Riemannian metric on $E$ is of bounded geometry (in $U$). We call it {\bf of bounded geometry} if for every section $\g$ there is an $\e > 0$ such that for the $\e$-neighborhood $U_{\e}$ of $\g (M) \subset E$, the bundle $\pi$ is of strongly bounded geometry in $U_{\e}$. 
\end{Definition}

\v Every bundle with compact fibers and compact base is of strongly bounded geometry, every bundle with compact base is of bounded geometry.

\v We can reformulate the conditions in the case of vector bundles:

\begin{Theorem}
A vector bundle $\pi: E \rightarrow M$ is of bounded geometry if and only if the injectivity radius of the base is bounded from zero, if the curvature of the base satisfies $\vert (\nabla^{(M)}) ^{(k)} R^M \vert \leq C^{(M)}_k$ for all $k \in \N \cup \{ 0\}$, and if the vector bundle curvature satisfies $\vert (\nabla^{(\pi)}) ^{(k)} R^{( \pi ) } \vert \leq C_k$ for all $k \in \N \cup \{ 0\}$. 
\end{Theorem}

\v {\bf Proof.} First we treat the curvature conditions: The tensoriality of $\nabla^{(k)}$ and of $R$ allows us to choose all horizontal fields basic and all vertical fields translation-invariant. Then we use the theorem \ref{haha}, \ref{huhu} and \ref{hoho} to show that the bounds of $\vert \nabla ^{(k)} R \vert$ in the definition above is equivalent to the existence of such bounds for the curvature terms of the corresponding vector bundle connections respectively bounds of the corresponding terms in the base. The geodesics split in two cases. Either $\dot{c} (0) $ is vertical; then $c$ stays in the fibre and $\dot{c} (t) = \dot{c} (0)$ modulo translation. Or it has a horizontal part in which case the geodesic lies over the geodesic $k: \In \rightarrow M$ given by $k(0) = \pi (c(0))$ and $\dot{k} (0) = d \pi (\dot{c} (0))$, and one can easily check by means of Theorems \ref{huhu} and \ref{haha} that the geodesic equations for $c$ are equivalent to the ordinary differential system $\nabla_k \iota (\dot{c}^{vert})=0$ where $\iota$ is the usual isomorphism between the vertical space at a section and the total space of a bundle. One only needs to cover the image of $k$ by convex neighborhoods and extend the vertical part of $\dot{c}$ translational-invariant. Thus $c$ is defined as long as $k$ is defined. Therefore around $c(0)$ there is a cylindrical neighborhood $\pi^{-1} (B_R( \pi (c(0))))$ which contains a ball of radius $R$ around $c(0)$.   \hfill \qed

\bigskip

\v There are some useful embedding theorems for sections of bundles of bounded geometry. For example, the compact supported smooth sections $\Gamma_0^{\infty} (\pi) $ are dense in the space of $L_p$-sections, for $1 \leq p < \infty$; and for $ s>k + \frac{n}{p}$ there is an embedding of $W^{s,p}$ sections into the space $C^{k}_b (\pi)$ of $k$ times differentiable sections of finite $C^k$-norm, cf. \cite{a}. 

\bigskip

\v Another useful property of bounded geometry is the following beautiful theorem which ensures the existence of a bounded atlas:

\begin{Theorem}
\label{nk}
Let $(M,g)$ be a Riemannian manifold of dimension $n$ which is of bounded geometry. Then there is a ball $B$ around $0$ in the Euclidean space $\R^n$ such that 

(i) $exp_p \vert_B$ is a diffeomorphism at every $p \in M$,

(ii) The distortion $\d (exp_p \vert_B)$ is a bounded function on $M$,

(iii) There is a bound $K_l$ independent of $p \in M$ such that all Christoffel symbols $\Gamma_{ij}^k  $ with respect to the normal neighborhood around all points $p \in M$ have $C^i$-norm smaller than $K_l$ on $B \subset \R^n$.    
\end{Theorem}

\v For the {\bf proof} cf. \cite{roe}, Lemma (2.2) and Proposition (2.4) as well as the proof in \cite{abp}, Appendix 2.  \hfill \qed

\begin{Theorem}
\label{bounded2}
Let $(M,g)$ and $(N,h)$ be two Riemannian manifolds and $F:M \rightarrow N$ be a smooth map. Then for every compact subset $K \subset M$ there are constants $\k_{K,i}, \hat{\k}_{K,i}$ such that

$$ \langle (\nabla^N)^{i-1} (F_* V), (\nabla^N)^{i-1} (F_* V) \rangle_{(\tau^*_N)^i \x \tau_N}  \leq  \hat{\kappa}_{K,i} \cdot  \langle (\nabla^M)^i  V, (\nabla^M)^i  V \rangle_{(\tau^*_M)^{(k)} \x \tau_N}  $$

\v and

\bean
 \langle \nabla^N_{F_* X_1}  \nabla^{N}_{F_* X_2} ... \nabla ^N _{F_* X_i} (F_* V ), \nabla^N_{F_* X_1}  \nabla^{N}_{F_* X_2} ... \nabla ^N _{F_* X_i} (F_* V )  \rangle_{ \tau_N}\\  
\leq  \kappa_{K,i} \cdot  \langle \nabla^M_{X_1} \nabla^M_{X_2} ... \nabla ^M _{X_i} V , \nabla^M_{X_1}   \nabla^M_{X_2}   ... \nabla ^M _{X_i}  V   \rangle_{ \tau_M} 
\eean

\v for any vector fields $V, X_1...X_i$ on $K$ where the covariant derivatives are the Levi-Civita ones.
 \end{Theorem}

\v {\bf Proof.} First we take a covering of every compact region $K_n$ by coordinate neighborhoods as in Theorem (\ref{nk}). Then we pick arbitrary coordinates $(x_1...x_{m-1}, t= x_{m})$ at a given point $(p, \tau) \in M \times \In$. Then we show the second claim first: Show inductively that $ \nabla^M_{X_1} \nabla^M_{X_2} ... \nabla ^M _{X_k} V$ is a sum of at most $k! \cdot 4^k \cdot m^{2k+1}$ terms of the form 

$$f^1_{i, j_1^1...j^1_{q(1)}} \cdot ... \cdot  f^n_{i, j_1^n...j^n_{q(n)}} \cdot \Gamma^{k_1}_{i_1 j_1, J^1_1 ... J^1_{Q(1)}} \cdot ... \cdot \Gamma^{k_p}_{i_p j_p , J^p_{1} ... J^p_{Q(p)}} V_{j, k_1...k_r} \partial_{x_n} $$

\v (where $X_i = \sum_j f^j_i \partial_j$ with $q(i), Q(i), p, r  \leq n$; this can be done by the observation that each factor in front of $\partial_{x_n}$ can be derived producing $m$ terms of one factor more and $\partial_{x_n}$ can be derived producing $m^2$ terms of two factors more). Then consider the pushed-forward coordinate system and apply the coordinate formula for the Christoffel symbols in partial derivatives of the metric. The partial derivatives of the coefficients $f_i^j$ stay the same on both sides. 

For the first part of the claim note that 

$$ \nabla^{n+1} (X_1, ... X_{n+1}) V = \nabla_{X_{n+1}} (\nabla^{n} (X_1, ... X_{n}) V) - \sum_{i=1}^{n} \nabla^{n} (X_1, ...\nabla_{X_{n+1}} X_i ... X_{n}) V $$ 

\v and show, again by induction, that $\nabla^{n} (X_1, ... X_{n}) V$ is a sum of $n!$ terms of the form

$$\nabla_{\nabla_{V(1,1)}... \nabla_{V(I(1), 1)} V(1)  } ...\nabla_{\nabla_{V(1,m)}... \nabla_{V(I(m), m)} V(m)  } V  $$

where $V(i,j), V(l) \in \{ X_1,... X_n \}  $ and if $k$ denotes the index we have 

$k( V (i,j)) > k(V(i+1,j)) > k(V(j)) > k(V(j+1)) $. Thus $\hat{ \k }_k  := ( n!  \max_{i \in 0,... k} \k_i )^2$ does the job. \hfill \qed

\section{Spaces of smooth sections and Fr\'echet geometry}
\label{frechet}

\v In this section we provide the foundations of Fr\'echet spaces and Fr\'echet manifolds, thus the material presented is novel only when dealing with concrete metrics. The non-metric aspects might be found in \cite{rh} or in \cite{cb}, however, we will include some of their proofs here.

\subsection{Fr\'echet spaces}

\v A motivation for generalizing the notion of Banach spaces is the fact that there is no Banach space structure on the space of smooth sections of a vector bundle such that the covariant derivative in the direction of any nonzero vector field on the base be continuous. This is seen by noting that in Banach spaces continuity is equivalent to boundedness, and by constructing sections which are arbitrarily high eigenvalues of the derivative or its square (most easily seen in $C^{\infty} (\S^1, \R)$ with $f_K (x) = K^{-1} sin (K^2 x)$).

\begin{Definition}
\label{basic}
A {\bf Fr\'echet space} is a locally convex topological vector space $V$ whose topology can be induced by a complete translational-invariant metric $d$ on $V$. If we keep track of its metric we call the pair $(V,d)$ {\bf metric Fr\'echet space}. For a real number $K$, a Fr\'echet space is called {\bf scalar-bounded by $K$} iff $d( \rho \cdot v, 0) \leq K \rho d(v, 0)$ for every $\rho \geq 1$.
\end{Definition}

\v {\bf Remark.} The triangle inequality implies that every Fr\'echet space with star-shaped balls is scalar-bounded by 2. However, even in finite-dimensional metric vector spaces, balls do not have to be star-shaped. As an example, consider the real line with the metric $d(r,s) := \Phi (\vert r-s \vert)  $ with $\Phi (x) := x$ for $0 \leq x \leq 1$, $\Phi (x) := 1- (x-1)/2$ for $1 \leq x \leq 2$ and $\Phi (x) := 1/2 + (x-2)/3$ for $x \geq 2$. as $\Phi (x \pm y) \leq \Phi (x) + \Phi (y)$, the metric $d$ satisfies the triangle inequality, but the balls with radius $1/2 \leq r \leq 1 $ are not starshaped and not even connected in this example.

\bigskip

\v {\bf Example 1:} Every Banach space $(V, \vert \vert \cdot \vert \vert)$ is a metric Fr\'echet space scalar-bounded by $1$, e.g. finite-dimensional vector spaces, $W^{k,p}$ and $C^p$, with $d(v,w):=\vert \vert v-w \vert \vert$.

\bigskip

\v {\bf Example 2:} The vector space $F^{\N}$ of sequences in a fixed metric Fr\'echet space $(F,d)$ can be made a metric Fr\'echet space $(F^{\N}, D)$ scalar-bounded by $1$ by setting 

$$\d_n((v_i)_{i \in \N}, (w_i)_{i \in \N}) := \sum_{i=1}^n d(v_i, w_i) $$

\v picking, for a $D>0$, a concave monotonously increasing continuous function $\Phi:[0, \infty )$ with $\Phi(0) = 0$ (from now on we will take the special choice $\Phi: [0, \infty) \rightarrow [0, 1), \Phi (x) = x/(1-x) $.), a positive sequence $\a$ converging to $0$, and defining the {\bf standard $ \a$-metric} 

$$D_{\a} ((v_i)_{i \in \N}, (w_i)_{i \in \N}) := \sum_{n \in \N} \a_n \Phi (  \d_n ((v_i)_{i \in \N}, (w_i)_{i \in \N}))$$

\v as well as the {\bf supremum $ \a$-metric}

$$d_{\a} ((v_i)_{i \in \N}, (w_i)_{i \in \N}) := sup_{n \in \N} \a_n \cdot \Phi(\d_n((v_i)_{i \in \N}, (w_i)_{i \in \N}) ) .$$

\v Each of these choices defines a natural family of metrics for the adequate category of Riemannian fiber bundles with connection. The supremum metrics have the advantage to have convex balls while they are nowhere differentiable. In contrast, as $\Phi$ is differentiable on $ [ 0, \infty)$, the standard metrics have better differentiability properties but none of their balls is convex. An interesting open question is whether there is a metric with both advantages. In contrast to Gl\"ockner (\cite{hg}) we choose the name {\em standard metric} for the {\em first} one as the usual metric appearing in the literature is $D_{\a}$ for $\a_n := 2^n$. For $r>0$ let $l(r)$ be the sequence with $l(r)_n = r^n$, then we have the estimate

\bea
\label{comparable}
D_{l(r^2)} < d_{l(r)} < D_{l(r)}  
\eea

\v which allows us to translate results obtained in the supremum metrics into the context of standard metrics and vice versa. The second one of these inequalities is obvious, for the first one note that as $l(r)$ tends to zero and the range of $\Phi$ is bounded, the supremum $s$ is attained by, say, $r^n \Phi (a_n)$. Then for all $j \in \N$ we have $\Phi (a_j) \leq r^{-n+j} s $ and therefore

$$
D_{l(r^2)}(a) = \sum_{i=1}^{\infty} r^{-2i} \Phi (a_i) \leq \sum_{i=1}^{\infty} r^{-2i} r^{-n+i} s \leq r^{-n} s \sum_{i=1}^{\infty} r^{-i} = r^{-n} s  .
$$

\v Note, however, that the inverse estimates $D_{l(r)} \leq C d_{l(r)}$ or $d_{l(r)} \leq C D_{l(r^2)}$ that would complete (\ref{comparable}) to equivalences between metrics, are wrong.
 
\bigskip

\v {\bf Example 3:}

\v Let $(M,g)$ be a Riemannian manifold and $\pi : E \rightarrow M$ a vector bundle over $(M,g)$ equipped with a Riemannian vector bundle metric $\langle \cdot , \cdot \rangle$ and a metric covariant derivative. Now we define metrics on $\Gamma_p (\pi)$ by pull-back of the above standard and Gl\"ockner metrics applied on the series of $L^p$ norms in the case of a compact base (by convention, put $L^{\infty} = C^0$):

$$\vert \vert \g \vert \vert^{(p,n)}  := \vert \vert \sqrt{ \langle \nabla^{(n)} \g ,  \nabla^{(n)} \g \rangle_{ (\tau ^* _M )^{\x n} \x \pi} } \vert \vert_{L^p (M)} $$

$$ \vert \vert \g \vert \vert _{p,n}: = \sum_{i=1}^n \vert \vert \g \vert \vert ^{(p,i)} $$

\bea
\label{metric}
\langle \g \rangle_{p,\a} &: =& D_{\a} ((\vert \vert \g \vert \vert _{p,n})_{n \in \N}, 0),\\
\langle \g \rangle^s_{p, \a} &: =& d_{ \a} ((\vert \vert \g \vert \vert _{p,n})_{n \in \N}, 0)
\eea

\v giving rise to $D_{p, \a}$ and $d_{p, \a}$, respectively, by applying $\langle \g \rangle_{p, \a}$ and $\langle \g \rangle^s_{p, \a}$ to differences. Likewise, if $M$ is noncompact, we define in a similar way standard and supremum metrics on $\Gamma(\pi) $ using an increasing series $C_n \subset C_{n+1}$ of compacta whose union is $M$. As we assume $M$ to be Riemannian and complete, we restrict ourselves to the choice $C_n := \overline{B}_n(x) $ where $x $ is an arbitrarily fixed point in $M$. Then we define

\bea
\label{metric2}
d_{x, p, \a }(\g, \d) := \sum_{i=1}^{\infty} 2^{-i} d_{x, p, \a}^i (\g, \d) &,& \qquad  d_{x, p, \a}^i (\g, \d) := d_{\a} \{ \vert (\g - \d) \vert_{K_j} \vert \vert_{W^{i,p}} \}_{j \in \N}, 0)  ,\\
D_{x, p, \a }(\g, \d) := \sum_{i=1}^{\infty} 2^{-i} D_{x, p, \a}^i (\g, \d) &,& \qquad  D_{x, p, \a}^i (\g, \d) := D_{\a} \{ \vert (\g - \d) \vert_{K_j} \vert \vert_{W^{i,p}} \}_{j \in \N} , 0)  .
\eea

\v The estimates (\ref{comparable}) immediately carry over to this case. 

\v {\bf Remark.} We could have made Definition (\ref{metric}) instead of (\ref{metric2}) also in the case of $M$ noncompact, if we assume that $\pi$ bounded and if we restrict ourselves to the spaces of sections all of whose Sobolev norms are finite. These are metric Fr\'echet spaces as well, natural in the class of (bundles over) Riemannian manifolds and isometries {\em without} base points; but as the restriction to these subspaces has some disadvantages we stick to (\ref{metric2}) in the noncompact case.

\begin{Theorem}
Let $\pi:E \rightarrow M$ be a vector bundle, then the space of its smooth sections $\Gamma_p (\pi)$ with the metrics $d_{p, \a}$ or $D_{p, \a }$ (in the compact case) resp. $d_{x, p, \a}$ or $D_{x, p, \a}$ is a metric Fr\'echet space scalar-bounded by $1$. The resulting topology is the same for all these metrics and is finer than the compact-open topology, it equals the inverse limit topology of the Banach spaces of $W^{k,p}$ sections in the compact case. The supremum metrics have convex balls, while none of the balls of the standard metrics is convex.
\end{Theorem}

\v {\bf Proof.}  It is enough to show that the metrics (i) are translation-invariant, (ii) are complete and has (iii) the balls in the supremum metrics are convex (then the topology is locally convex in either case). 

\bigskip

\v Property (i) is clear by definition. For (iii), we have to show that for each two sections $\g, \g '$ we have $\langle \g + t (\g ' - \g) \rangle^s_{p, \a} \leq max\{ \langle \g \rangle^s_{p, \a}, \langle \g ' \rangle^s_{p, \a} \}$ for all $ t \in [0,1] $. By continuity of the metric it is sufficient to restrict ourselves to the case $t = 1/2$, and the rest is then done by nesting by intervals. By concavity of the arithmetic sum and monotonicity of $\Phi$ we have 

$$\Phi ( \frac{a_n + b_n}{2}) \leq max \{ \Phi (  a_n), \Phi (  b_n)  \}   $$

\v and therefore, for $d^s_{\a, p}$, 

$$\a_n \Phi (  \frac{a_n + b_n}{2}) \leq max  \{  \Phi (  a_n), \a_n \Phi (  b_n)  \}     $$

\v and thus 

$$sup(\a_n \Phi ( \frac{a_n + b_n}{2})) \leq max  \{ sup (\a_n \Phi (  a_n)) , sup( \a_n \Phi (  b_n))  \}     .$$

\v For Property (ii), consider a Cauchy sequence of smooth sections. Then we have to show that there is a limit in $\Gamma (\pi)$. But

$$ \langle \d \rangle \geq \sum_{n=0}^k \Phi (\a_n \vert \vert \d \vert \vert_{n,p}) \cdot 2^{-n} \geq 2^{-k} \sum_{n=0}^k \Phi (\a_0 \vert \vert \d \vert \vert_{n,p}) $$

\v Therefore a Cauchy sequence $\g_n$ in $\Gamma (\pi)$ is a Cauchy sequence in the Banach space $W^{k,2}$ as well (put $\d = \g_k - \g_l$ and use the continuity of $\Phi$ at $0$ and the invertibility of $\Phi$ in the positive semiaxis) and has therefore a limit $\g$ in this Banach space. Because of uniqueness this limit is the same in all these Banach spaces. As we can use the Sobolev embedding theorems in every $C_n$, it is smooth. For the statement about the compact-open topology note that the form of the standard metric shows that a set $A$ is open in the Fr\'echet space if it is open in every $\Gamma^k_p (\pi)$. The metrics $D_{p, \a}$ and $D_{x, p, \a}$ are scalar-bounded by $1$ as every of its additive terms is. The metrics $d_{p, \a}$ and $d_{x, p, \a}$ are scalar-bounded by $1$ as well because if for fixed $v \in \G (\pi) $ we consider the real function $L: t \mapsto d_{p, \a} (tv, 0)$. As this function consists piecewise of concave functions which vanish in $0$, it is easy to see that $L(s \cdot t) \leq s L(t)$.   \hfill  \qed

\bigskip

\v In all the examples the metric could be constructed by a countable family of seminorms. This is a general feature of Fr\'echet spaces as shown by the following definition and the following theorem that can be found in \cite{rh}, \cite{cb}:

\begin{Definition}
A {\bf seminorm} $\vert \vert \cdot \vert \vert $ on a vector space $F$ is a real-valued, nonnegative function which is subadditive and satisfies $\vert \vert cf \vert \vert = \vert c \vert \cdot \vert \vert f \vert \vert$ for all scalars $c$ and vectors $f$. 
\end{Definition}

\begin{Theorem}[cf. \cite{rh}, e·g.]
Let $F$ be a Fr\'echet space. Then there is a $\N$-family of continuous seminorms $\vert \vert \cdot \vert \vert_i$ on $F$ whose balls $B_{\e}^i (x) := \{ y \in X: \vert \vert y-x \vert \vert _i < \e \}$ are a basis of the topology of $F$. Therefore the topology of $F$ can be generated by the metric 

\bea
\label{sum}
D_{\a} (f,g) := \sum_{i=1}^{\infty} \a_n \Phi ( \vert \vert f-g \vert \vert_i)
\eea

\v where $\a$ is an arbitrary positive sequence converging to $0$, and as well by

\bea
\label{sup}
d_{\a} (f,g) := \sup_{i \in \N} \a_n \Phi ( \vert \vert f-g \vert \vert_i).
\eea

\end{Theorem}

\v {\bf Proof.} Choose a Fr\'echet metric $d$, consider $B_{\frac{1}{i}}^d (0) $ and define the seminorms as the so-called Minkowski functionals 

$$ \vert \vert v \vert \vert _i := inf \{ \l > 0 \vert \frac{1}{\l} \cdot v \in  U_i  \} .$$

\v where we choose convex subsets $U_i \subset B_{\frac{1}{i}}^d (0)$. These Minkowski functionals are subadditive, as for $\frac{1}{\l} f, \frac{1}{\mu} g \in U_i$ we have also $\frac{1}{\l + \mu} (f+g) \in U_i$ as a convex combination. Continuity is an easy consequence of subadditivity. Finally, Cauchy sequences w.r.t. all $\vert \vert \cdot \vert \vert _i$ are Cauchy sequences for the metric.    \hfill \qed

\bigskip

\v It should be stated, however, that if the metric $d$ was already given as a sum of seminorms as in Equation (\ref{sum}), the Minkowski functionals will {\em not} give us back the original seminorms (although the tame equivalence class stays the same as we will see in Theorem \ref{minkeq}).
One reason why the description by seminorms is important is that it appears in the Nash-Moser inverse function theorem which is valid only in the case of {\em tame} Fr\'echet spaces and {\em tame} operators which in turn are defined in terms of the seminorms. Thus let us recall some definitions from the tame category:

\begin{Definition}
A linear map $A$ between Fr\'echet spaces $F,G$ with sequences of seminorms  $\vert \vert \cdot \vert \vert_i $ and $\vert \vert \cdot \vert \vert_i ' $, respectively, is called {\bf tame (from $(F,\vert \vert \cdot \vert \vert_i )$ to $(G,\vert \vert \cdot \vert \vert_i ')$)} if there are natural numbers $b, r$ such that $A$ satisfies estimates $ \vert \vert v \vert \vert_n ' < C_n \vert \vert v \vert \vert_{n+r} $ for $n>b$.  
Two different sequences of seminorms $\vert \vert \cdot \vert \vert $ and $\vert \vert \cdot \vert \vert ' $ on one and the same Fr\'echet space $F$ are called {\bf tamely equivalent} if there is a linear isomorphism $A:F \rightarrow F$ which is tame from $(F, \vert \vert  \cdot \vert \vert_i)$ to $(F, \vert \vert \cdot \vert \vert_i' )$ and whose inverse is tame from $(F, \vert \vert  \cdot \vert \vert_i')$ to $(F, \vert \vert \cdot \vert \vert_i )$.   
\end{Definition}

\begin{Theorem}
\label{minkeq}
Let $F$ be a Fr\'echet space with a metric defined by Equation \ref{sum} or Equation \ref{sup}. Then the sequence of norms defining the metric is tamely equivalent to the sequence of $2^i$-Minkowski functionals.
\end{Theorem}

\v This theorem (whose straightforward proof shall be skipped) may sound quite technical, but it implies that tameness of maps does not depend on the choice of a special sequence of seminorms but only on the metric as the latter one gives rise to the sequence of Minkowski functionals. Thus tameness is an inherent notion of the {\em metric} category and will be related to the sequence of $2^i-$Minkowski functionals from now on. It even only depends on the equivalence class of a metric which is shown by the following theorem whose proof will be omitted as well:

\begin{Theorem}
Let $A: (F,d) \rightarrow (G,d')$ be bounded by $C$, then $A$ is tame with grade $r \in \N$ whenever $2^r >C$, base $0$ and $C_n =1$ for all $n \in \N$.
\end{Theorem}

\v Let us come back to our examples. Comparing Example 2 and Example 3 we notice that in Example 2 none of the seminorms we used is a norm while in Example 3 any of the seminorms is a norm. The question could arise whether there is {\em any} continuous norm on $F^{\N}$. This question is answered negatively in the following theorem.

\begin{Theorem}[cf. \cite{rh}]
The Fr\'echet space $F^{\N}$ does not have a continuous norm.
\end{Theorem}

\v {\bf Proof.} Let us assume the existence of a continuous norm $\nu$. Then we consider a ball $B_R^{\nu} (0)$. On one hand, this ball cannot contain any nontrivial subspace of $F^{\N}$, as $\nu$ is homogeneous w.r.t. the multiplication by positive numbers. But on the other hand, the ball is open because of continuity of $\nu$, so it contains a finite intersection of elements $B_R^{\vert \vert \cdot \vert \vert _i} (0)$ of the basis of the topology. But an intersection of the balls for the seminorms $\vert \vert \cdot \vert \vert _{i_1}, ... \vert \vert \cdot \vert \vert _{i_n}$ contains the subspace $\{ x \in F^{\N} \vert x_1 = ... = x_m =0 \} $ where $ m = max \{ i_1, ... i_n \}$, a contradiction. \hfill \qed

\bigskip

\v {\bf Example 4:} Spaces of maps ($C^k$, $W^{k,p}$ or smoothly finite) from a Riemannian manifold to a metric Fr\'echet space $F$ can be made metric Fr\'echet spaces (with the same scalar bound). This example can be generalized to sections of Fr\'echet vector bundles.

\bigskip

\v A theorem which provides us with still more examples for metric Fr\'echet spaces (and whose non-metric variant can be found in \cite{rh} or \cite{cb}) is:

\begin{Theorem}
(1) A closed subspace of a metric Fr\'echet space is again a metric Fr\'echet space, scalar-bounded by the same constant. 

(2) A quotient of of a metric Fr\'echet space by a closed subspace is again a metric Fr\'echet space, scalar-bounded by the same constant. 

(3) The direct sum of finitely many metric Fr\'echet spaces is again a metric Fr\'echet space, scalar-bounded by the maximum of the bounds.

\end{Theorem}

\v {\bf Proof.} (i) Restrict the metric to the subspace and consider the relative topology of the closed subspace. Convex sets stay convex as intersected with a linear subspace. The scalar bound is trivial.

\v (ii) Let us call the closed subspace $U$ and the surrounding Fr\'echet space $X$. Define the new metric $d'$ by $d' (v,w) := min_{c \in U} d( v +c, w) =  min_{c,d \in U} d( v +c, w+d)$ (the last equation is valid because of the invariance of $d$ under translations). This metric generates the quotient topology. Now for every Cauchy sequence in $X/U$ we have to find a Cauchy sequence of representatives in $X$. Thus choose a $M_{\e} \in \N$ s.t. for all $m,n >M$ we have $d' ([v_m], [v_n]) = min_{c \in U} d(v_m, v_n +c) < \frac{\e}{3}$. Then choose a $m(0) > M$, a representative $v_{m(o)}$ and a sequence of vectors $c_n \in U$ with $d (v_{m(0)}, v_n + c_n) < \frac{\e}{2}$. Then using the triangle inequality we see that for $\tilde{v}_n := v_n + c_n$ we have $d(\tilde{v}_k, \tilde{v}_l) < \e$. Now modify the sequence of representatives successively this way for $\e = \frac{1}{n}$ for all $n \in \N$. this converges and leaves us with a Cauchy sequence in $X$. For the scalar bound and for $\rho \geq 1$ take $\tilde{c} := \rho \cdot c $ in the definition of the distance.

\v (iii) Let $d_1, d_2$ the two metrics, we choose a continuous concave function $\Delta:\R^2 \rightarrow \R $ (e.g. $x_1 + x_2$ or $\sqrt{x_1^2 + x_2^2}$) and define the new metric $d' := \Delta \o (d_1, d_2)$. For the scalar bound use concavity of $\Delta$.  \hfill \qed

\bigskip

\begin{Theorem}[metric Hahn-Banach-theorem]
Let $F$ be a Fr\'echet space, $G \subset F$ a subspace and $\l : G \rightarrow \R $ a continuous linear map. Then there is a continuation of $\l$ to a continuous linear map $F \rightarrow \R$. If we fix a Fr\'echet metric $d$ with respect to which $\l$ is bounded on $G$ by $R$, we can choose a continuation bounded by $R$ as well. In particular, for every vector $f \in F$, there is a continuous linear functional $\l$ on $F$ with $\l (f) \neq 0$.
\end{Theorem}

\v {\bf Proof.} The proof is in complete analogy to the Banach case, cf. \cite{dw}, pp. 372-378. For the question of boundedness take the Theorem of Hahn-Banach as quoted in cf. \cite{dw}, pp.94-97: Let $X$ be a real vector space, $U$ a subvector space of $X$, $p:X \rightarrow \R$ sublinear, $l:U \rightarrow \R$ linear with $l(u) \leq p(u)$, for all $u \in U$. Then there is a linear continuation $L:X \rightarrow \R$, $L \vert_U = l$, with $L(x) \leq p(x)$, for all $x \in X$. To apply this, take $p(x) = sup_{u \in U} \frac{\vert \vert l(u) \vert \vert}{d(u,0)} \cdot d(x,0) $.  \hfill \qed
 
\bigskip

\begin{Theorem}[Open-mapping-theorem]
Let $F,G$ be Fr\'echet spaces. If $T:F \rightarrow G$ is a linear map which is continuous and surjective, then $T$ is an open map. In particular, if $T$ is continuous and bijective, $T$ is a homeomorphism. 
\end{Theorem}

\v {\bf Proof} analogously to the Banach case, cf. \cite{ky}, p.75 \hfill \qed

\bigskip

\begin{Definition}
A subspace $G$ of a Fr\'echet space $F$ is called {\bf topologically complemented in $F$} if there is another subspace $H$ of $F$ such that the map $\Phi: G \times H \rightarrow F, \Phi ((g,h)) := g + h$ is a homeomorphism, i.e. if $F$ is homeomorphic to the topological direct sum $G \oplus H$. In this case we call $H$ a {\bf topological complement of $G$ in $F$}. 
\end{Definition}

\v As there is a homeomorphism from $G$ to $G \oplus \{ 0 \} \subset G \oplus H$, a topological complemented subspace is always closed. Likewise, it is easy to see that $G$ is topologically complemented in $F$ if and only if there is a continuous projection $\pi : F \rightarrow G$.

\bigskip

\v {\bf Example:} Hamilton (\cite{rh}) gives an example of a closed subspace of a Fr\'echet space which is not topologically complemented. So take $F:= C^{\infty} ([0, 1])$ which contains the space $G$  of $1$-periodic real functions on the real line, $C_1^{\infty} (\R) $ by the restriction $\rho$ on the unit interval. If we define 

$$p: C^{\infty} ([0,1]) \rightarrow \R^{\N}, f \mapsto D^j f (2 \pi ) - D^j f (0)$$

we get the short exact sequence

$$\{ 0 \} \rightarrow C^{\infty}_1 (\R) \rightarrow^{\rho}  C^{\infty} ([0,1]) \rightarrow^p \R^{\N} \rightarrow \{ 0 \}$$

\v Thus the quotient of $C^{\infty} ([0,1]) $ by $ C^{\infty}_1 (\R)$ is homeomorphic to $\R^{\N}$. As the latter one does not have any continuous norm, there cannot be a continuous linear isomorphism between $\R^{\N}$ and any closed subspace of $F$. Therefore the above sequence does not split, and $G$ is not topologically complemented in $F$.

\bigskip

\v This behaviour is not exceptional which is shown by the following theorem:

\begin{Theorem}[cf. \cite{gk}, p. 435]
Let $F$ be a Fr\'echet space with a continuous norm which is not Banach. Then there is a closed subspace $H \subset F$ with $F/H \cong \R^{\N}$, thus $H$ is not topologically complemented in $F$. 
\end{Theorem}

\v But at least simple subspaces of Fr\'echet spaces are topologically complemented:

\begin{Theorem}[\cite{gk}]
\label{topcom}
Let $F$ be a Fr\'echet space. Then

(1) Every finite-dimensional subspace of $F$ is closed.

(2) Every closed subspace $G \subset F$ with $codim (G) = dim (F/G) < \infty$ is topologically complemented in $F$.

(3) Every finite-dimensional subspace of $F$ is topologically complemented.

(4) Every linear isomorphism between the direct sum of two closed subspaces and $F$, $G \oplus H \rightarrow F$, is a homeomorphism.       
\end{Theorem}

\subsection{Differentiation and integration of Fr\'echet maps}

\v Contrary to the case of Banach spaces we will have several different notions of differentiability of maps between Fr\'echet spaces. First define the spaces $L(F,G)$ of set-theoretic linear maps and $CL(F,G)$ of continuous linear maps  between Fr\'echet spaces $F,G$. Then we note that $CL(F,G)$ can be made a topological vector space by the compact-open topology. 

\begin{Theorem}
Let $F, G$ be topological vector spaces, let $G$ be locally convex and metrizable. Then the set $CL(F,G)$ of continuous linear maps topologized by the compact-open topology is a topological vector space. 
\end{Theorem}

\v {\bf Proof.} Choose a metric $d$ inducing the topology of $G$. Let $A+B=C$ in $CL(F,G)$, let $C \subset (K,O)$. Then for $M:= A(K), N:= B(K)$ we know that $M,N$ and therefore also $M+N$ compact and $M+N \subset O$. Now let $\e:= d(M+N, \partial O) > 0$. Then $B_{\frac{\e}{2}} (M) + B_{\frac{\e}{2}} (N) \subset O $ by the triangle inequality, and $(K, B_{\frac{\e}{2}} (M)) \times (K, B_{\frac{\e}{2}} (N))$ is an open neighborhood of $(A,B)$ which is mapped in $(K,O)$ under $+$. Thus $+$ is continuous. Now let us prove that the scalar multiplication is continuous as well. Let $r \in \R, A \in CL(F,G)$ be given as well as an open neighborhood of $rA$ of the form $(K,O)$. We have to find a positive number $s$ and an open neighborhood $U$ of $A$ with $(r-s, r+s) \times U$ mapped into $(K,O)$ by scalar multiplication. First note that there is an $\e >0$ with $B_{\e} (rA(K)) \subset O$. Now consider the continuous function $R: B \mapsto \sup_{x \in K} d(A(x), B(x))  $ on $CL(F,G)$ and define $U:= R^{-1}( [ 0, \frac{\e}{8r}   ))$. Choose a symmetric convex open neighborhood $W$ of $0$ contained in $B_{\e / 2} (0)$ and consider the continuous Minkowski functional $\rho (v) := inf \{ t \in \R \vert \pm t^{-1} v \in W \} $ on the compact set $AK$ where it attains a finite maximum which we call $s'$. Put $s:= min \{ s', r \} $. Now $(r-s, r+s) \times U$ has the desired property: Let $(\tau, B) \in (r-s, r+s) \times U $, let $k \in K$, then

\bean
d(\tau Bk, rAk) &\leq& d(\tau Bk, \tau Ak) + d(\tau Ak, rAk)\\
&\leq& S \tau d(Bk, Ak) + d((\tau- r) Ak, 0)\\
&<& \frac{\e}{2} + \frac{\e}{2} = \e  
\eean

\v where $S$ is a bound of scalar multiplication which can be chosen smaller or equal to $2$. This concludes the proof. \hfill \qed 

\bigskip

\v At a first glance, this looks good: Why not define differentiability between Fr\'echet spaces by means of the spaces $CL(F,G)$? The overnext theorem will tell us that this concept would not be very far-reaching as it does not allow for the iterative definition of higher differentiability. As a preparation for its proof we define certain subsets of $\G (\pi)$ which will turn out to be compact and generic for compact subsets in the sense of the next theorem. For $\g \in \Gamma(\pi) $ and an $a \in \R^{\N}$ we set 

$ K^{\g, \a} = \{ s \in \G (\pi) \vert \ \vert \vert \b - \g \vert \vert_n \leq a_n  \ {\rm for \ all}\  \b \in B \ {\rm and \ all} \ n \in \N$.

\v As an intersection of closed sets, every $ K^{\g, \a}$ is closed. The following theorem shows that it is compact as well.

\begin{Theorem}
A closed set $B \subset \Gamma(\pi)$ is compact if and only if there is a $\g \in \Gamma(\pi) $ and an $a \in \R^{\N}$ with $B \subset K^{\g, \a}$.
\end{Theorem}

\v {\bf Proof.} Assume the condition fails, then there is an $n \in \N$ such that $B$ is unbounded w.r.t. $\vert \vert \cdot \vert \vert_n$ and, consequently, contains a sequence $\b_n$ which in turn does not contain any Cauchy subsequences in this norm. Therefore none of its subsequences converges in $\vert \vert \cdot \vert \vert_n$ and thus neither in the Fr\'echet space $\G (\pi)$. On the other hand, it is easy to check that each $K^{\g, \a}$ is closed in $\G (\pi)$. Now, for a given $K^{\g, \a}$, define $K_l^{\g, \a} := \{ \d \in \Gamma^{l+1} \vert \ \vert \vert \d- \g \vert \vert_i \leq \a_i \forall i = 1, ... l+1  \} $. The closure $\overline{K_l^{\g, \a}}$ of $K_l^{\g, \a}$ in $\G^l (\pi)$ is compact because of the Arzela-Ascoli Theorem. The series $\overline{K^{\g, \a}}$ converges to $\overline{K^{\g, \a}}= K^{\g, \a}$ in the inverse limit $\G (\pi) = \lim_{\leftarrow} \G^l (\pi)$, thus by the compatibility of compactness of Hausdorff spaces and inverse limits the latter one is compact in $\Gamma(\pi)$ (this fact can be looked up in \cite{ems} or \cite{du}, App.2, 2.4., or seen by the fact that the inverse limit is defined as a subset of the Tychonoff product which is compact in this case and it is a closed subset by the constituents of the limit being Hausdorff). \hfill \qed

\bigskip

\begin{Theorem}
$CL(F,G)$ is not metrizable.
\end{Theorem}

\v {\bf Proof.} Recall that every metric space is first-countable, i.e. every point has a countable neighborhood base. Now take $0 \in CL(F,G)$ and assume that $0$ has a countable neighborhood base $U_n$. Every $U_n$ is the union of finite intersections of sets of the form $(K_i, O_i)$ where the $O_i$ are open neighborhoods of $0 \in G$. So for every $n$ we pick one finite intersection of the union and put $\tilde{U}_n := (\bigcup K_i, \bigcap O_i)$. Note that $\tilde{U}_n \neq \empty $ as $\bigcap O_i$ is an open neighborhood of $0$. If $\{ U_n \}_{n \in \N} $ was a neighborhood base so is $ \{ \tilde{U}_n \}_{n \in \N} $ as $\tilde{U}_n \subset U_n$ for every $n$. Now put $U_n ' := (K^{\g_n, \a_n}, \bigcap O_i )$ where $K^{\g_i, \a_i}$ contains $\bigcup K_i$. W.r.o.g. we can choose $\g_i = 0$ for all $i$. Then, following Cantor's diagonal procedure, define $\omega_k:= 2 (\a_k)_k $ and $K:= K^{\g, \omega}$. This is a neighborhood of $0$ which contains none of the system above. This concludes the proof by contradiction. \hfill \qed

\begin{Definition}
Let $F,G$ be two Fr\'echet spaces, $U \in F$ open. A map $Q : U \rightarrow G$ is called {\bf differentiable at $p \in U$} iff there is a linear map $A_p: F \rightarrow G$ with 

$$ A_p(v) = \lim_{t \rightarrow 0} \frac{Q(p + tv) - Q(p)}{t} $$

\v for all $v \in F$. If $Q$ is differentiable at all points $p \in U$, if $dQ(P)$ is continuous for all $p \in U$ and if the induced map $dQ: U \rightarrow CL(F,G)$ is continuous, we call $Q$ {\bf c-differentiable}. The space of all c-differentiable maps from $U$ to $G$ is denoted by $C^1(U,G)$.
\end{Definition}

\begin{Theorem}
For an open set $U \in F$ a map $Q :U \rightarrow G $ is c-differentiable iff it is differentiable at every point and if the map $Q' : U \times F \rightarrow G , (u,f) \mapsto dQ (u) f$ is continuous in the product topology.
\end{Theorem}

\v {\bf Proof.} Continuity of $df$ at every point is trivial. Let $df: U \rightarrow CL(F,G)$ be given, and without restriction of generality let $df(p) =0$. Let an element of the basis $V=(K,O) \ni 0 $ be given. We have to show that there is an open neighborhood $U$ of $p$ with $df(U) \subset V$ or equivalently $f'(u, K) \subset O$ for all $u \in U$. By continuity of $f'$, around every point $q \in K$ we can find an open neighborhood $W_q$ of $q$ and an open neighborhood $U_q$ of $p$ with $f' (U_q \times W_q) \subset O$. Let $W_1, ... W_n$ be a finite subcovering of $K$, then for $U: = \bigcap_{i=1}^n U_i$ we have

$$ f'(U \times K) = f'(\bigcap_{i=1}^n U_i \times K) \subset f' (\bigcap_{i=1}^n U_i \times \bigcup_{i=1}^n  W_q) \subset O,$$

\v thus $U$ has the required property. The other direction is easily seen as well e.g. in the treatment of the exponential law in \cite{du}.   \hfill \qed

\bigskip

\v The last result allows us to extend the notion of c-differentiability to the notion of {\bf Keller differentiability} or {\bf k-differentiability} for short, which {\em will} allow for iterative definitions of higher derivatives:

\begin{Definition}
Let $F,G$ be two Fr\'echet spaces, $U \in F$ open. If $Q: U \rightarrow G$ is differentiable at all points $p \in U$, if $dQ(P)$ is continuous for all $p \in U$ and if the induced map $Q': U \times F \rightarrow G$ is continuous, we call $Q$ {\bf k-differentiable}. The space of all k-differentiable maps from $U$ to $G$ is denoted by $K(U,G)$. Inductively we can define $Q^{(k+1)}: U \times F^{k+1} \rightarrow G$ by $Q^{(k+1)} (u, f_1, ... f_{k+1}) = \frac{d}{dt} Q^{(k)} (u + t f_{k+1} , f_1, ... f_k  $, for every $k$ as this consists of continuous maps and is continuous. The space of maps $Q: U \rightarrow G$ for which $Q^{(i)}$ exists is denoted by $K^i (U, G)$. 
\end{Definition}

\begin{Definition}
Let $(F,d), (G,d')$ be two metric Fr\'echet spaces, $r \in (0, \infty]$. The space $B_r(F,d,G,d')$ is the space of all $r$-bounded linear maps between $F$ and $G$, that means, all maps $f$ for which

$$ d^{(r)}_{F,G} (f,0) :=  \langle f \rangle_r := \sup_{p \in B_r(0) \setminus \{ 0\} } \frac{d(f(p), 0)}{d(p,0)}  $$

\v is finite. If no confusion can occur, we will denote this space also by $B_r(F,G)$.
\end{Definition}

\v {\bf Remark.} Of course always $r \leq R \Rightarrow B_R (F,G) \subset B_r(F,G)$. If $d_G$ is finite, for every $\e >0$, any given map is bounded on $F \setminus B_{\e} (0)$ by $\e^{-1}$, thus $B_r(F,G) = B_R(F,G)$ as subsets of $CL(F,G)$ for all $r, R \in \R \cup \{ \infty \} $, as for $A \in B_r(F,G)$ we have $\langle A \rangle_R \leq max \{ \langle A \rangle_r, r^{-1} \}  $. Nevertheless, although for $d_G$ finite they coincide as subsets of $CL(F,G)$, the spaces $B_r(F,G)$ carry different topologies. Note that $B_r(F,G)$ can be made a (locally convex, if $(G, d')$ has arbitrarily small convex balls) complete metric vector group with metric $\d (f, g):= \langle f-g \rangle$; it is, however, in general {\em not} a Fr\'echet space and not even a topological vector space as was pointed out in \cite{hg} where it was proven that the obvious (pointwise) scalar multiplication is not continuous in $B_{\infty}(F,G)$ at $t=0$ if $F$ or $G$ contains a line through $0$ on which the metric is bounded. Actually one can formulate an even more saddening result:

\begin{Proposition}
Let $F$, $G$ be Fr\'echet spaces. Let $d_1$ be a metric on $F$ and $d_2$ be a metric on $G$ which generate the respective topologies. Then for every $r\in (0, \infty] $, if $G$ is not Banachable and $B_r(F,d_1,G, d_2)$ is a topological vector space, it does not contain any surjective map with a bounded right inverse. If  $F$ is not Banachable and $B_r(F,d_1,G, d_2)$ is a topological vector space, it does not contain any map bounded away from zero. In the particular case of $(F,d_1) = (G, d_2)$, $B_r(F,d,F,d)$ is not a topological vector space, unless $F$ is Banachable.  
\end{Proposition}

\v {\bf Proof.} First we show that for a non-Banachable metric Fr\'echet space $H$, for every real numbers $K,r$ there are vectors $h \in B_r (0)$ with $l_K(h) := d(Kh)/d(h) $ arbitrarily close to one. Namely, if we could bound $l_K$ from $1$ by, say, $\e >0$, take $U_0$ a convex reflection-symmetric open neighborhood of $0$ contained in $B_1 (0)$ and (let w.r.o.g. $K>1$) $U_n := \{ h \in H \vert K^n \cdot h \in U_0 \} $. As $U_n \subset B_{(1-\e)^n} (0)$, this would be a neighborhood base for $0$, and it is easy to check that the norm $\vert \vert v \vert \vert := inf \{ \l \in \R \vert \l^{-1} \cdot v \in U_0 \} $ would generate the topology of $H$. Therefore we cannot bound $l_K $ from $1$, and there are $v_n \in B_r(0)$ with $l_K(v_n) < 1+ \frac{1}{n} $. Now consider first the case that $B_r(F,G)$ contains a surjective map $f$ with bounded right inverse $g$. Then $C:= im(g)$ is a Fr\'echet subspace of $F$ and $f \vert_{C}$ is an isomorphism bounded both from $\infty$ and $0$. But $\l f$ does not converge to $0$ as $\l \rightarrow 0$, because 

\bean
sup_{x \in B_r(0) \setminus \{ 0 \}} \frac{d( \l f(x), 0)}{d(x,0)} &=& sup_{x \in B_r(0) \setminus \{ 0 \}} (\frac{d_2(\l f(x) , 0)}{d_2(f(x), 0) } \cdot \frac{d_2 (f(x), 0)}{d_1(x, 0)})\\
&\geq& sup_{x \in B_r(0) \setminus \{ 0 \}} (\frac{d_2(\l f(x) , 0)}{d_2(f(x), 0) }) \cdot inf_{x \in B_r(0) \setminus \{ 0 \}}(\frac{d_2 (f(x), 0)}{d_1(x, 0)})
\eean

\v in which expression the first factor is $1$ and the second one independent of $\l$. So $B_r(F,G)$ is not a topological vector space in this case. 

\v In the other case  observe

\bean
sup_{x \in B_r(0) \setminus \{ 0 \}} \frac{d( \l f_2(x), 0)}{d_1(x,0)} &=& sup_{x \in B_r(0) \setminus \{ 0 \}} \frac{d_2(f(x), 0)}{d_1(\l^{-1} x, 0)}\\
&=& sup_{x \in B_r(0) \setminus \{ 0 \}} (\frac{d_2(f(x), 0)}{d_1(x, 0)} \cdot \frac{d_1(x, 0)}{d_1(\l^{-1} x, 0)} )\\
&\leq& sup_{x \in B_r(0) \setminus \{ 0 \}} (\frac{d_2(f(x), 0)}{d_1(x, 0)}) \cdot inf_{x \in B_r(0) \setminus \{ 0 \}} (\frac{d_1(x, 0)}{d_1(\l^{-1} x, 0)} )
\eean

\v and then proceed as above. For the statement about the case $F=G$ note that the identity is surjective and bounded from zero, independently of the metric used. \hfill \qed

\bigskip

\v Now, as in general $B_r(F,G)$ is not a topological vector space, again we face the problem how to define higher derivatives. $B_r(F,G)$ is an algebraic subspace of $L(F,G)$. But as its topology is finer than the one induced from $L(F,G)$ it is a nontrivial question whether it is open or closed in $L(F,G)$. It turns out that in general it is never open nor closed: Take $F= \R^{\N}$ with a standard metric, and for all $j \in \N$ consider $f_j: F \rightarrow F$ with $f_j (v_i) = i v_i$ for all $i \leq j$ and zero on the complement $C_j$ of the canonical embedded $\R^j \subset \R^{\N}$. Clearly in $L(F,F)$ the sequence converges to $f$ with $f (v_i) = i v_i$ for all $i \in \N$ and zero on $C$. But all $f_j$ are contained in $B_r(F,F)$ while $f$ is not. Thus $B_r(F,F)$ is not closed in $L(F,F)$. On the other hand, $0$ is contained in $B_r(F,F)$, and $F_j:= f- f_j$ converges to $0$ in $L(F,F)$ while all $F_j$ are unbounded. Therefore $B(F,F)$ is not open in $L(F,F)$, either.

\bigskip

\v We recall that, following Gl\"ockner's result, for Fr\'echet spaces $F,G$, the space $B_r(F,G)$ is not a topological vector space in general, but a complete metric vector group with absolutely convex balls (if we consider real Fr\'echet spaces, this means that the balls are invariant under reflections $x \mapsto -x$. In the complex case, they have to be invariant under multiplication with a complex unit). We call such a vector group a {\bf strong vector group}. It is easy to see that if $F, H$ are strong vector groups then $B_r(F, H)$ is again a strong vector group. Although we cannot really use $B_r(F,G)$ in order to {\em define} higher derivatives as this requires a topological scalar product, we can use the Keller derivatives and {\em reinterpret} them in the bounded sense: Define inductively strong vector groups $B^{(i+1)}_r(F,G):= B_r(F, B^{(i)}_r(F,G))$ and spaces of set-theoretic linear maps $L^{i+1}(F,G):= L (F, L^i(F,G))$. Then, for $PL(U \times F^k, G)$ denoting the space of maps which are linear in every argument except the first one, consider natural maps $T^k: PL (U \times F^k, G) \rightarrow L^k (F, G) ^U $. This enables us to make the following definition:

\begin{Definition}
Let $F,G$ be Fr\'echet spaces, let $U \subset F$ be open, and let $Q: U \mapsto G$ be a map. $Q$ is called {\bf bounded-differentiable} or {\bf b-differentiable with radius $R$} if it is differentiable at every point, if all maps $dQ(p)$, $p \in U$, are bounded within the ball around $0$ of radius $R$ and if the map $dQ: U \rightarrow B_R(F,G)$ is continuous.  For $i \in \N \cup \{ \infty \}$, we define $Q$ to be $i$ times b-differentiable with radius $r$ in $U$ if $Q \in K^i (U, G) $ we have that all $T^k (Q^{(k)}) : U \rightarrow L^k (F,G)$ take their image in $B_r^{(k)} (F,G)$ for $k \leq i$ and are continuous w.r.t. the strong vector group topology. The space of all $i$ times b-differentiable maps with radius $R$ we denote by $ B_R^i (U,G)$.  
\end{Definition}

\v Of course the notion of b-differentiability is the strongest one and implies k-differentiability in very much the same way as uniform continuity implies continuity. A trivial corollary from the definition is

\begin{Theorem}
Let $A: F \rightarrow G$ be an affine map, $A(f) = L f + g $ for all $f \in F$ where $L$ is a linear map and $g \in G$. Then $A \in K^{\infty} (F,G)$ and $A'_f = L$ for all $f \in F$. Consequently, the map $A$ is b-differentiable with radius $r$ if and only if $L$ is $r$-bounded.  \hfill \qed
\end{Theorem}

\v H.Gl\"ockner proved in \cite{hg} that the topology of $B_{\infty}(F,G)$ is totally disconnected. The connectivity properties of $B_r(F,G)$ for finite $r$ is a little bit better:

\begin{Theorem}
The connected component $B_{r,0} (F,G)$ of $0$ in $B_r(F,G)$ contains all maps in $B_r (F,G)$ which map $B_r(0)$ into a compact set. If $G= \G (\pi)$ with a metric of the form \ref{sum} or \ref{sup}, $B_{r,0} (F,G)$ equals the set of  all maps in $B_r (F,G)$ which map $B_r(0)$ into a compact set. 
\end{Theorem}

\v {\bf Proof.} From \cite{hg}, Prop. 2.10, we recall that $A \in B_{r,0} (F,G) $ if and only if $\lim_{t \rightarrow 0} d(tA, 0) =0$. If $A (B_r(0))$ is precompact, all Minkowski functionals are bounded on $A(B_r(0))$, thus it is easy to see that the limit is $0$. If $G= \G (\pi)$, then there is a $C^k$ norm which is unbounded on $A(B_r(0))$. Therefore the corresponding Minkowski functional is unbounded as well, so the limit is not zero. \hfill \qed

\begin{Theorem} 
\label{fine}
For a vector $v \in \G (\pi)$, the curve $c:t \mapsto t \cdot v $ is b-differentiable if and only if there is a real number $M$ such that $\vert \vert v \vert \vert_k \leq M $ for all $k \in \N$.
\end{Theorem}

\v {\bf Proof.} The derivative of $c$ is $\dot{c}(t) (s \partial_t):= s \cdot v$, and one finds 

$$\sup \langle sv \rangle/s = \lim_{s \rightarrow 0} \langle sv \rangle / s = d/ds \vert_{s= 0} \langle sv  \rangle.$$  

\v The corresponding sum converges exactly if the sequence of seminorms is bounded. \hfill \qed

\bigskip

\v As in \cite{rh} explained, for any Fr\'echet space $F$, the Riemannian integral along curves can be defined by its universal property that it commute with linear functionals; it has similar proprties as in finite-dimensional or Banach analysis, e.g. the fundamental theorem of calculus carries over to the Fr\'echet case. Likewise one can establish many of the usual theorems like the chain rule for differentials. For more on this, see \cite{rh}.

\v H.~Gl\"ockner proved a useful estimate linking the distances $d_F$, $d_G$ and $d_{F,G}$. The assumptions of this theorem contains the requirement that $d_G$ have absolutely convex balls. 

\begin{Theorem}[\cite{hg}, Lemma 1.11]
\label{convexbound}
Let $(F, d_F)$ and $(G,d_G)$ be metric Fr\'echet spaces such that $d_G$ has absolutely convex balls. Let $U \subset B_R (x_0) \subset F$ be a convex subset with non-empty interior and $f \in C^1( U,  G)$. Then for all $x,y \in U$ we have

$$d_G(f(x), f(y)) \leq d_F(x,y) \cdot  \sup_{t \in [0,1]} \langle f' (x+t(y-x)) \rangle_R .$$ 
\end{Theorem}

\v Now, an easy consequence of the established terminology which can be proved in complete analogy to the finite-dimensional case with the usual $3 \e$-arguments is the following theorem:

\begin{Theorem}
\label{convergence}
Let $F$, $G$ be Fr\'echet spaces, let $U \subset F$ be open, let $f_n: F \rightarrow G$ be a pointwise convergent sequence of $B^1_r (U,G)$ maps with one and the same bound $B$ and let all $df_n$ be Lipschitz functions with one and the same Lipschitz constant $B$. Then the map $f$ of pointwise limits is bounded-differentiable in $U$ with the same bound $B$, and $lim_{n \rightarrow \infty} f'_n = f'$. \hfill \qed
\end{Theorem}

\v {\bf Remark.} If all $f_n$ are $B_r^2 (U,G)$ with the same bound $B$, the Lipschitz condition is satisfied automatically.

\bigskip

\v Now we want to show that there are interesting examples of smooth bounded maps between Fr\'echet spaces:

\bigskip

\v {\bf Example 1:} 

\begin{Theorem}
Let $\pi: E \rightarrow M$ be a vector bundle with a metric and a metric connection. For a parallel vector field $X$ on $M$, the map $\nabla_X: \Gamma(\pi) \rightarrow \Gamma (\pi)$ is bounded by $2 \vert \vert X \vert \vert_0 $ in both standard and supremum metrics. 

\end{Theorem}

\v {\bf Proof.} As $\langle \nabla^{(k)} (\nabla_X \g) \rangle \leq \langle \nabla^{(k+1)} \g \rangle$, the problem reduces to the boundedness of the shift $\s, ( \s(a) )_n := a_{n+1}$ in $\R^n$ which is easily shown to be bounded by $2$. \qed

\begin{Theorem}
For every unbounded metric generating the standard Fr\'echet topology on $\R^{\N}$ the shift is unbounded as well. 
\end{Theorem}

\v {\bf Proof.} Let $d$ be the metric in question. The sequence of onen subsets $U_n := \{ a \in \R^{\N} \vert a_i < 1/n \ {\rm for} \  {\rm all} \ i \leq n \} $ is a neighborhood basis of $0$. Therefore there is an $m \in \N $ with $U_m \in B_1^d (0)$. Choose $v \in \R^{\N}$ with $d(v, 0) > \s^m$. But 

$$(\frac{1}{2m}, \frac{1}{2m}, ... \frac{1}{2m}, v_1, v_2, ... ) \in U_m \subset B_1^d(0) $$

\v (where the first $m$ coordinates are meant to be $1/2m$) which is a contradiction. \hfill \qed

\bigskip

\v {\bf Example 2:} Let $\pi, \psi$ be vector bundles over $M$ and let $A \in \Gamma ((\pi^*)^n \times \psi)$ (a polynomial with values in $\psi$). If $A$ is parallel w.r.t. the canonical connection in the tensor product, then for every collection of parallel vector fields $X_i^j$ on $M$ the map

$$C_A: \g \mapsto A(\nabla_{X_1^1} ... \nabla_{X_1^{n(1)}} \g , ... ,\nabla_{X_m^1 } \nabla_{X_m^{n(m)}} \g )    $$

\v is $B_r^{\infty}$ for every radius $r$. This is easily seen from the facts that polynomials are bounded in every ball, that the covariant derivative is bounded as above, and that the contraction with a parallel tensor $A$ is bounded as $ \nabla_X (A( \g)) = A(\nabla_X \g)    $, therefore $ \nabla^{n} A( \g) = A(\nabla^n \g  )  $, and $ \langle C_A \rangle \leq \vert \vert A \vert \vert_0 $. One example for this construction is $\g \mapsto \langle \gamma, \nabla_X \g \rangle  $ for a parallel vector field $X$.

\bigskip

\v {\bf Counterexample:} The most striking counterexample for boundedness is the composition of maps. Let $f: \R^n \rightarrow \R^m$ a smooth map, let $M$ be a compact manifold and let $C_f: C^{\infty}(M, \R^n ) \rightarrow C^{\infty} (M, \R^m)$ be the composition with $f$. It is well-known that $C_f$ is a smooth tame map, but even in the case of $f$ having compact support it is in general wrong that $C_f$ is b-differentiable in standard or supremum metrics. This can easily be seen by taking $ f:= \Phi \cdot sin (M x_1)$ where $\Phi$ is a compactly supported function which is identically $1$ in a neighborhood of the origin.

\subsection{Fr\'echet manifolds}

\v Due to its better differentiability we will perform the following steps only for standard metrics although they can in principle also be done for Gl\"ockner metrics.

\begin{Definition}
A Hausdorff topological space $\mathcal{M}$ is called a {\bf Fr\'echet manifold} if $\mathcal{M}$ has an open covering $U_{\a}, \a \in A$ such that for every $\a \in A$ there exists a homeomorphism $\Phi_{\a}: U_{\a} \rightarrow V_a $ to an open subset $V_a$ of a Fr\'echet space $F_{\a}$ such that for any $\a, \b \in A$ with $U_{\a} \cap U_{\b} \neq \emptyset$ the map

$$ \Phi_{\a} \o \Phi_{\b} \vert_{\Phi_{\b} (U_{\a} \cap U_{\b})}: \Phi_{\b} (U_{\a} \cap U_{\b}) \rightarrow \Phi_{\a} (U_{\a} \cap U_{\b})  $$

\v is a smooth map between two open sets of Fr\'echet spaces. If there is a Fr\'echet metric in each of these Fr\'echet spaces and if we require additionally the maps $ \Phi_{\a} \o \Phi_{\b} \vert_{\Phi_{\b} (U_{\a} \cap U_{\b})}$ to be bounded w.r.t. these metrics, we will speak of a {\bf bounded Fr\'echet manifold}. If there is a uniform bound for these chart transitions we speak of a {\bf strongly bounded Fr\'echet manifold}.  \end{Definition}

\begin{Definition}
A {\bf compatible metric} on a Fr\'echet manifold $M$ is a metric $d$ on $M$ such that there is a Fr\'echet subatlas of $M$ such that in each chart $U$, $d$ is equivalent to the Fr\'echet metric $d_U$, that is, there are constants $b_U, B_U$ with $b_U \cdot d_U(p,q) \leq d(p,q) \leq B_U \cdot d_U (p,q)$ for any two points in $U$, and if $d$ is scalar-bounded in every chart. If there is a choice of $b_U, B_U$ independent of $U$ we speak of a {\bf strong compatible metric}.
\end{Definition}

\v This definition implies that the metric $d$ generates the topology of $M$. Moreover, an important property of $d$ implied by the definition is that for every point of $M$ there is a real number $R_p$ such that Cauchy sequences in a ball $B_r (p)$ have a limit in the closure of the ball as long as $r \leq R_p$. An easy consequence of the definition of a compatible metric is the following theorem:

\begin{Theorem}
If a Fr\'echet manifold carries a (strong) compatible metric, then it is (strongly) bounded.
\end{Theorem}

\v To make use of the notion of Fr\'echet manifolds (which we want to model also the spaces of sections of bundles over non-compact manifolds), we have to make sure that the bundles we are working with are geometrically not too wild. To this purpose we use the notion of bounded geometry introduced in section \ref{bounded}.

\begin{Theorem}
Let $(M,g)$ be a Riemannian manifold and $\pi : E \rightarrow M$ a fiber bundle over $(M,g)$ equipped with a Riemannian fiber bundle metric $\langle \cdot , \cdot \rangle$ and a fibre bundle connection $D$ such that the bundle is of (strongly) bounded geometry. Then the space of smooth sections $\Gamma (\pi)$ can be given the structure of a bounded Fr\'echet manifold with a (strong) compatible metric.
\end{Theorem}

\v {\bf Proof.} First we define a metric on $\Gamma(\pi)$ whose values are bounded by $1$.
For two sections $\g, \d \in \Gamma (\pi)$ which are homotopic to each other by a smooth homotopy $H \in \Gamma(\pi \times [0,T])$ of sections $\g_t= H(\cdot, t)$ of $\pi$, we define

$$H_m (p,t) := \langle  \nabla^{(m)} \partial_t H (p,t),  \nabla ^{(m)} \partial_t H (p,t) \rangle_{(\tau^* _M)^{\x m} \x \tau (\pi) }$$

\v in which formula $\nabla = \nabla ^{\g_t^* \tau_E}$ which is well-defined as all $\g_t$ are immersions into $E$ being sections of $\pi$.

\v Having defined $H_m (p,t)$, we define

$$ H ^{(n,k)} := \sum_{m =0}^n  \vert \vert \sqrt{  \vert \vert H_m \vert \vert_{L^k (M_t)} }  \vert \vert _{L^1([0,T])} $$

$$\langle H \rangle^{(k, \a)} := \langle \{ H^{(i,k)} \}_{i \in \N} \rangle_{\a}  $$

\v and finally

$$d^{(k, \a)}(\g_0, \g_1) :=  \inf_{H: \g \leadsto \d} \{ \langle H^{(k)} \rangle \} $$

\v We will omit the indices $(k, \a)$ whenever this does not cause any confusion. Sometimes, if we refer to a special metric $g$ in the fibre, we will use the notion $d_g$. Obviously, the function $d^{(k)}$ is {\bf symmetric} by definition. It is {\bf positive}, as all its defining terms are positive and as the function $\Phi: x \mapsto \frac{x}{1+x}$ maps the positive semiaxis to itself. The vanishing of $d^{(k)}$ implies that $\g = \d$ as $H^{(0,k)} \geq \int dist_f(\g (m) , \d (m)) dm $ where $dist_f$ denotes the fibrewise Riemannian distance. This quantity is greater than zero for $\g \neq\d$. 
The {\bf triangle inequality} follows by puzzling together two isometries: if two sections $\g_0, \g_1$ are homotopic to each other by an homotopy $H_{01}$ then we can find another homotopy $H'_{01}$ between them all of whose $t$-derivatives vanish at $t= 1$ (by concatenating $H_{01}$ with a diffeomorphism $\psi_{01}$ of $[ 0,1 ]$ all of whose derivatives vanish at $1$) and still with $\langle H'_{01} \rangle_p = \langle H_{01} \rangle_p$ by reparametrization-invariance of the arc-length quantities $H^{(n,k)} $. The same works for an homotopy $H_{12}$ between $\g_1$ and $\g_2$; here the reparametrized homotopy $H'_{12}$ having vanishing $t$-derivatives at $t=1$. Now defining an homotopy $H_{02}$ between $\g_0$ and $\g_2$ by the prescription $H_{02} := H'_{01}$ for $t \in [0,1]$ and  $H_{02} := H'_{12}$ for $t \in [1,2]$, we can show by elementary calculus that $H_{02}^{(n)} = H_{01}^{(n)} + H_{12}^{(n)}$ for all $n$, and the convexity of $\Phi : x \mapsto \frac{x}{1+x}$ implies $\langle H_{02} \rangle \leq \langle H_{01} \rangle + \langle H_{12} \rangle$ which in turn by infimum arguments implies the triangle inequality.

\v The {\bf Fr\'echet manifold structure} is provided by smoothly finite sections (w.r.t.$ d_{\g}$) of $\g ^* \tau_M^v$ around a section $\g$ of $\pi$. The topology will be induced by $U \in \tilde{\Gamma}(\pi)$ open if $\kappa^{-1} (U)$ open for all charts $\kappa$ of the atlas.

\v We equip an open set of the pulled-back vertical bundle $\g ^* \tau^v_E$ with two different metrics $g_{\g}$ resp. $g$ by pullback along $\g$ resp. along the identification with an open set in $\Gamma(\pi)$ by the exponential map in the fibre. Then we have to prove the following proposition:

\begin{Proposition}
The metric $d$ and the translational-invariant pullback metric $d_{g_{\g}} =: d_{(\g)}$ are topologically equivalent as metrics on sufficiently small $d_{(\g)}$-open subsets $U$ of $\tilde{\Gamma}_{d_{(\g)}} (\g ^* \tau_M^v)$, i.e. every $d_{(\g)}$-ball around a section in $U$ contains a $d$-ball and vice versa. 
\end{Proposition}

\v {\bf Proof.} First we need to prove a lemma which shows in the same time that in the case of $\pi$ being a vector bundle with translational-invariant metric and vector bundle connection the metric we define is the same as the metric defined in the section about Fr\'echet spaces:

\begin{Lemma}
\label{liniso}
If we equip $\g ^* \tau_M^v$ with a vector bundle metric $h$ (as e.g. the pullback metric $g_{\g}$) and a metric vector bundle connection (as e.g. the pullback connection) we have $d_{g_{\g}}(\a, \b) = \langle L \rangle_h$ where $L=(1-t) \a + t \b$ is the affine homotopy joining $\a$ and $\b$.   
\end{Lemma}

\v {\bf Proof.} As by definition of a vector bundle metric every fibre is a flat vector space with a translational-invariant metric, the minimizing homotopy has to be affine as one can see by decomposing the homotopy in an affine part and one perpendicular to it: First we put w.r.o.g. $p=0$. As $l$ is a continuous convex function it is enough to show that the straight line segment is a {\em local} minimum of $l$, therefore we can restrict ourselves to the open space of curves $c: 0 \leadsto q$ with $\langle q , \dot{c} \rangle >0 $. Let $\e>0$ be given. Then pick $I,N \in \N$ such that $2^{-I} 2^{-N} < \e $. Consider the decomposition $\dot{c} = f q + s $ in the Hilbert space determined by the norms up to $I,N$ where $s \perp q $ in this Hilbert space (and so in all scalar products taking part in it). As the length functional is invariant under reparametrizations, we can parametrize these curves $c$ such that $f = 1 $. Then we have $ \langle c \rangle \geq \langle g \rangle - \e = d(0, q) - \e $.  \hfill \qed

\bigskip

\v {\bf Proof of the proposition, continued.} It is enough to show that for all sections $\a \in \tilde{\Gamma}_{g_{\g}} (\g ^* \tau_M^v)$ and for all $\d >0$, there is a $\rho >0$ such that $B^{d_{(\g)}}_{\rho} (\a) \subset B^{d}_{\d} (\g) $ and vice versa, i.e. with the roles of $d$ and $d_{(\g)}$ exchanged. First we show the inequalities between the metrics: For the charts we use Riemannian normal coordinates restricted on the open subsets appearing in Theorem \ref{nk}.

\bigskip

\begin{Lemma}
\label{X}
Let $\pi$ be of bounded geometry, pick one section $\g$. Then there are constants $\k_n$ with $I_{g_{\g}}^{i,k} \leq \k_i I_g^{i,k} $ and $I_{g}^{i,k} \leq \k_i I_{g_{\g}}^{i,k} $ for all isotopies $I$ which satisfy $dist^h (\g, I) \leq \b$ where $\b$ is the radius of the ball in Theorem \ref{nk}. 
\end{Lemma}

\v {\bf Proof of the lemma.} First show that the parallel transport along fiber geodesics is uniformally bounded using Theorem \ref{nk}. Then note that $I^* \tau^v_E$ and $\gamma^* \tau^v_E$ are exactly related by this parallel transport (which corresponds to radial translation in the exponential plane) and apply Theorem \ref{huhu} and Theorem \ref{bounded2}. \hfill \qed

\bigskip

\v Now we have the inequalities (for $L$ being the linear isometry in the chart linking $\a$ with $\b$, with $\k_i$ as above, $\bar{\k}_n := max_{i \in \{ 1, ... n\}} \k_i$)

\bean
d(\a, \b) &\leq& \langle L \rangle_g\\
&\leq& \sum_{i=1}^n 2^{-i}\Phi (L^{(n)}_g) + 2^{-n} \ \rm{(rest} \  \rm{sum} )\\
&\leq& \sum_{i=1}^n 2^{-i}\Phi( \k_i L^{(i)}_{g_{\g}}) + 2^{-n} \ {\rm(monotonicity} \ \rm{ of } \  \Phi , \ {\rm  Lemma} \  \rm{\ref{X})}\\
&\leq& \bar{\k}_n  \sum_{i=1}^n 2^{-i} \Phi (L^{(i)}_{g_{\g}}) + 2^{-n} \ {\rm (scalar-boundedness} )\\
&\leq& \bar{\k}_n \sum_{i=1}^{\infty} 2^{-i} \Phi (L^{(i)}_{g_{\g}}) + 2^{-n}\\               
&=& \bar{\k}_n \langle L \rangle_{g_{\g}} + 2^{-n} \ {\rm (Lemma} \ \rm{ \ref{liniso})}\\
&= \bar{\k}_n d_{\g} (\a, \b) + 2^{-n}
\eean

\v Likewise, for $I$ being an isotopy between $\a$ and $\b$ with 

\begin{equation}
\label{vorauss}
\langle I \rangle_g < 2 d (\a, \b) 
\end{equation}

\v we get the inequality

\bean
d_{\g} (\a, \b) &\leq& \langle I \rangle_{g_{\g}}\\
&<&  \sum_{i=1}^n 2^{-i} \Phi (I^{(i)}_{g_{\g}}) + 2^{-n} \ {\rm (rest} \ \rm{sum)}\\
&\leq& \sum_{i=1}^n 2^{-i} \Phi ( \k_n I^{(i)}_g ) + 2^{-n}\ \rm{(monotonicity} \ \rm{of } \  \Phi, {\rm Lemma} \ \rm{ \ref{X})}\\ 
&\leq& \bar{\k}_n \sum_{i=1}^n 2^{-i} \Phi (I^{(i)}_g) + 2^{-n} \ {\rm (convexity} \  \rm{of} \ \Phi)\\
&\leq& \bar{\k}_n \sum_{i=1}^{\infty} 2^{-i} \Phi (I^{(i)}_g) + 2^{-n}\\        
&\leq& \bar{\k}_n \langle I \rangle_g + 2^{-n} \ {\rm (Definition} \  \rm{of} \ \langle I \rangle _g)\\
&<& 2 \bar{\k}_n d (\a, \b) + 2^{-n} \ {\rm (Condition} \  \rm{( \ref{vorauss} ))}
\eean

\v Thus, the inequalities $ d(\a, \b) < 2 \k_n  d_{\g} (\a, \b) + 2^{-n} $ and $d_{\g} (\a, \b) < 2 \k_n  d (\a, \b) + 2^{-n}$ hold true for every given natural number $n$. Thus let a ball $B_R^{d} (\g)$ be given, then choose $n$ with $R> 2 ^{-n+1}$, then elementary calculus shows that $B^{d_{\g}}_{2^{-n-1} \k_n^{-1}} (\a) \subset B_R^{d} (\a)$ which proves the proposition. \hfill \qed

\bigskip

\bigskip

\v Now the modelling vector space is {\bf locally convex} as the topologies induced by the translational-invariant metric $d_{\g}$ and the one induced by $g$ in any chart are equivalent, as shown above.

\v The metric $d$ is {\bf locally complete} in the sense specified above because of the inequalities established above. But the Fr\'echet manifold itself is in general not a complete metric space (see below).   \hfill \qed

\bigskip

\v Note that this construction would not be possible if in the definition of the notion of Fr\'echet manifold the existence of a countable basis of the topology were required. 

\v {\bf Remark.} We defined the metric $d$ for spaces $\Gamma(\pi)$ of sections of fiber bundles. What does happen if we want to deal with maps $f: M \rightarrow N$ between manifolds? Of course, we would represent them as sections of the trivial fiber bundle $p_1: M \times N \rightarrow M$. Now let us consider two sections $\g, \d$ of $\Gamma(\pi) $ and forget about their property being sections, only considering them as maps from $M$ to the total space $N$ of $\pi$. If we denote by $d^t$ the metric in the space of sections of the associated trivial bundle $M \times N$, it is obvious that we get $d^t (\g, \d) \leq d(\g, \d)$ as every homotopy as a section can be seen as a homotopy as a map, and the metrics that apply coincide. The converse is wrong in general which can be seen in the example of $C^{\S^1, \R}$ by considering to sharp peaks centered at two points very close to each other but with disjoint supports : They are very close as maps but can be arbitrarily distant as sections.

\begin{Theorem}
Let $\pi: E \rightarrow M$ be a bounded bundle. If all fibers of $\pi$ are complete, then the Fr\'echet manifold $\Gamma(\pi)$ with the metric $d_g$ is a complete metric space. \end{Theorem}

\v {\bf Proof.} As $d_g$ contains the fiber distance as an additive term, if we evaluate a $d_g$-Cauchy sequence $\g_i$ at a point $p \in M$, we get a Cauchy sequence in $\pi^{-1} (p)$ which converges to a point $q$. By the usual estimates in uniform convergence one gets the continuity of the pointwise limit: By $C^1$-convergence of the $\g_i$ we can conclude the uniform continuity of the sequence and for a sequence $p_i \rightarrow p$ we can decompose $d_f(\gamma (p_i), \gamma (p) ) \leq d_f(\g (p_i) , \g_j (p_i)) + d_f(\g_j (p_i), \g_j (p) ) + d_f(\g_j(p), \g(p))$ to show that the pointwise limit $\g$ is continuous. Then consider the sections $\g_i$ of the Cauchy sequence now in a chart around $\g$ with the translational-invariant metric equivalent to $d_g$. Now for every vector field $V$ on $M$, the sequence of sections $\nabla_V \g_i$ is Cauchy and converges pointwise. Its continuity can be proven by $C^2$-convergence of $\g_i$. In this manner we can proceed inductively to prove the smoothness of $\g$.  \hfill \qed

\begin{Definition}
Let $(M,d), (N,d')$ be two Fr\'echet manifolds with compatible metric. The space $B_r(M,d,N,d')$ is the space of all $r$-bounded maps between $M$ and $N$, that means, all maps $f$ for which

$$\langle f \rangle_r := \sup_{p \neq q \in M, d(p,q) <r} \frac{d(f(p), f(q))}{d(p,q)}  $$

\v is finite. If no confusion can occur, we will denote this space also by $B(M,N)$.
\end{Definition}

\v This is a generalization of the same notion in the case of Fr\'echet spaces.

\begin{Theorem}
\label{ballisenough}
Let $R>0$ be given. In the case $(M,d)= (\Gamma(\pi), d_g), (N, d') = (\Gamma(\pi), d_{g'})$ and for a Fr\'echet map $f: M \rightarrow N$ we have

$$\langle f \rangle_R := \sup_{q \in M, p  \in B_r(q) \setminus \{ q \} } \frac{d(f(p), f(q))}{d(p,q)} .$$

\v for a $r <R$, i.e., we can restrict us to the case of pairs of points with any distance smaller than $R$. 
\end{Theorem}

\v {\bf Proof.} Let $\e >0$ be given. Define $\langle f \rangle_r := \sup_{q \in M, p \neq q \in B_r(q)} \frac{d(f(p), f(q))}{d(p,q)} $. Let $I$ be an isotopy between two points $A,B \in M$ with $\langle I \rangle < d(A,B) + \e$. Then consider the function $D(t) = \langle I \vert_{[0,t]} \rangle$ and split the isotopy into finitely many parts $I_i: a_i \leadsto a_{i+1}$ such that $\langle I_i \rangle < r $ for all $i$, in particular $d(a_i, a_{i+1}) < r$. Then we have $d'(f(a_i), f(a_{i+1})) <  \langle f \rangle _r d(a_i,a_{i+1})  $, and the iterative application of the triangle inequality implies that 

$$d' (f(A),f(B)) < \langle f \rangle_r \cdot \langle I \rangle < \langle f \rangle _r \cdot (d(A,B) + \e)  ,$$

\v for arbitrary $\e>0$, and the claim follows. \hfill \qed

\bigskip

\v Note that in the case of $M,N$ being Fr\'echet spaces with fixed Fr\'echet metrics and $f$ being linear, we have $\langle f \rangle_R  := \sup_{p \in B_r(0)} \frac{d(f(p), 0)}{d(p,0)} $ for all $r \leq R$.

\bigskip

\v Now let us rephrase the inverse function theorems of Nash and Moser in the language of the metric-tame category:

\begin{Theorem}[Left inverse function Theorem between Fr\'echet spaces]
\label{left}
Let $F$ and $G$ be metric-tame Fr\'echet spaces. Let $U \subset F$ be open, $f: U \rightarrow G$ smooth metric-tame, let $f' (x)$ be injective in $U$ with left inverse $L_0 :U \times G \rightarrow F $ which is metric-tame in $U$ as well. Then there is a real number $r_0 > 0$ such that there is a continuous left inverse of $f \vert_{B(x_0, r_0)}$ (thus in particular it is injective). \hfill \qed
\end{Theorem}

\begin{Theorem}[Right inverse function Theorem between Fr\'echet spaces]
Let $F$ and $G$ be metric-tame Fr\'echet spaces. Let $U \subset F$ be open, $f: U \rightarrow G$ smooth metric-tame with surjective differential in $U$ and smooth metric-tame right inverse $R_0: U\times G \rightarrow F$. Then there is a real number $r_0 > 0$ such that $V:= f(B(x_0, r_0))$ is open in $G$ and there is a smooth right inverse of $f \vert_{B(x_0, r_0)} \rightarrow V$ (thus in particular it is surjective). \hfill \qed
\end{Theorem}

\begin{Theorem}[Full inverse function Theorem between Fr\'echet spaces]
Let $F$ and $G$ be metric-tame Fr\'echet spaces. Let $U \subset F$ be open, $f: U \rightarrow G$ smooth metric-tame, let $f': U \times F \rightarrow G$ be a smooth metric-tame map which is an isomorphism at every $u \in U$. Then there is a real number $r_0 > 0$ such that $V:= f(B(x_0, r_0))$ is open in $G$ and $f \vert_{B(x_0, r_0)} \rightarrow V$ is a diffeomorphism. \hfill \qed
\end{Theorem}

\bigskip

\v As corollaries, we get the inverse function theorems for Fr\'echet manifolds with compatible metrics:

\begin{Theorem}[Left inverse function theorem between metric Fr\'echet manifolds]
Let $M$ resp. $N$ be Fr\'echet manifolds with compatible metrics $d_M$ resp. $d_N$. Let $x_0 \in U$, $U \subset M$ be open, $f: U \rightarrow N$ smooth metric-tame, let $T^1 (f') $ be injective in $U$ with left inverse $L: U \times G \rightarrow  F$ smooth metric-tame as well. Then there is a real number $r_0 > 0$ such that there is a continuous left inverse of $f \vert_{B(x_0, r_0)}$ (thus in particular it is injective). \hfill \qed
\end{Theorem}

\begin{Theorem}[Right inverse function Theorem between metric Fr\'echet manifolds]
Let $M$ resp. $N$ be Fr\'echet manifolds with compatible metrics $d_M$ resp. $d_N$. Let $x_0 \in U$, $U \subset M$ be open, $f: U \rightarrow G$ smooth metric-tame with surjective differential in $U$ and right inverse $R: U \times G \rightarrow F$ which is smooth metric-tame as well. Then there is a real number $r_0 > 0$ such that $V:= f(B(x_0, r_0))$ is open in $G$ and there is a smooth metric-tame right inverse of $f \vert_{B(x_0, r_0)} \rightarrow V$ (thus in particular it is surjective). \hfill \qed
\end{Theorem}

\begin{Theorem}[Full inverse function theorem between metric Fr\'echet manifolds]
Let $M$ resp. $N$ be Fr\'echet manifolds with compatible metrics $d_M$ resp. $d_N$. Let $U \subset M$ be an open set, $f: U \rightarrow N$ smooth metric-tame w.r.t. $d_M$ and $d_N$. Let $T^1 (f' ) $ be an isomorphism in $U$ with an inverse $I: U \times G \rightarrow F$ which is smooth metric-tame as well. Then there is a real number $r_0 > 0$ such that $V:= f(B(x_0, r_0))$ is open in $N$ and $f \vert_{B(x_0, r_0)} \rightarrow V$ is a diffeomorphism. \hfill \qed  
\end{Theorem}

\section{An inverse function theorem for bounded maps}

\v In this section we want to state an inverse function theorem in the category of Fr\'echet spaces and bounded maps. We will go very closely along the lines of a nice pedagogical introduction to the same subject in {\em Banach} spaces written by Ralph Howard (\cite{h}). The results are interesting due to the fact that unlike the Nash-Moser theorem invertibility is required only in one point. An interesting extension of the inverse function theorems given below was given by Helge Gl\"ockner in \cite{hg}, there only finite differentiability is required.

\begin{Theorem}[Von-Neumann series]
\label{von-neumann}
Let $(F,d)$ be a metric Fr\'echet space of finite metric and $A \in B_R(F, d, F, d)$ with $\langle A - \1_F \rangle_R = r < 1, r <R$. Then $A$ is invertible with inverse 

$$A^{-1} = \sum_{i=0}^{\infty} (\1_F - A)^i$$

\v which satisfies the estimate

$$ \langle A ^{-1} \rangle_R \leq \frac{1}{1 - \langle \1 _F - A \rangle_R }.$$

\v Moreover for $0 < \rho < 1/2$, $\langle A^{-1} \rangle < 1$, $\langle \1_F - A \rangle_R  , \langle \1_F - B \rangle_R \leq \rho$ and $\langle A - B \rangle_R < 1/2$ we have $\langle A^{-1} - B^{-1} \rangle_{R} \leq \frac{1}{(1- \rho)^2} \langle A - B \rangle_{2R}$.

\end{Theorem}

\v {\bf Proof.} Let $B:= \sum_{k=0}^{\infty} (\1_F - A ) ^k$. Then $\langle ( \1_F -A)^k \rangle_R \leq \langle \1_F - A \rangle_R ^k = r^k$, so the series defining $B$ converges by comparison with the geometric series on $B_R(0)$ and by linear continuation on all of $F$, and

$$\langle B \rangle_R  \leq \sum_{k= 0}^{\infty} \langle \1_F - A \rangle_R ^k = \frac{1}{1- \langle \1_F -A \rangle_R}  .$$ 

\v Now compute 

\bean
AB &= \sum_{k=0}^{\infty} A (\1_F - A)^k = \sum_{k=0}^{\infty} (\1_F - (\1_F -A)) (\1_f -A)^k\\
&= \sum_{k=0}^{\infty} (\1_F - A)^k - \sum_{k=0}^{\infty} (\1_F - A)^{k+1} = \1_F.
\eean

\v Likewise one can show $BA= \1_F$ or simply use that $A$ and $B$ commute. 

\v Now if $\langle \1_F - A \rangle_R, \langle \1_F - B \rangle_R \leq \rho $, the preceeding shows $\langle A^{-1} \rangle_R, \langle B^{-1} \rangle_R \leq \frac{1}{1- \rho} < 2$, and

\bean
\langle B^{-1} - A^{-1} \rangle_R &= \langle B^{-1} (B-A) A^{-1} \rangle_R \leq \langle B ^{-1} \rangle_R  \langle B-A \rangle_{R}  \langle A^{-1} \rangle_{R} \\
&\leq \frac{1}{(1 - \rho)^2} \langle A - B \rangle_{R}
\eean

\v which completes the proof. \hfill \qed

\begin{Theorem}
\label{neumann2}
Let $F$ and $G$ be metric Fr\'echet spaces and let $A,B \in B(F,G)$ and $A$ be invertible, $A^{-1} \in B_r(G,F)$ with $\langle A \rangle_r \leq 1$. Then if $\langle A -B \rangle_r < \frac{1}{2 \langle A^{-1} \rangle_r}, 1$, then $B$ is invertible as well and 

$$\langle B^{-1} \rangle_r \leq \frac{\langle A^{-1} \rangle_r}{1 - \langle A ^{-1} \rangle_r \langle A-B \rangle_r} $$

and

$$ \langle B^{-1} - A^{-1} \rangle_r \leq \frac{\langle A^{-1} \rangle_r ^2 \langle B -A \rangle_r }{1 - \langle A^{-1} \rangle_r \langle A-B \rangle_r } .$$

\v Thus the set $U$ of $r$-bounded-invertible maps $r$-bounded by $1$ from $F$ to $G$ is open in $\{ A \in B_r(F,G) \vert \langle A \rangle_r < 1 \}$, and the map $A \mapsto A^{-1} $ is continuous on $U$.
\end{Theorem}

\v {\bf Proof.} As $B= A(\1_F - A^{-1} (A-B))$ and $\langle A^{-1} (A-B) \rangle_r  \leq \langle A^{-1} \rangle_r \langle A-B \rangle_r < 1$ by assumption, Theorem \ref{von-neumann} shows that $\1_F - A^{-1} (A-B)$ is invertible and that

$$\langle (\1_F - A^{-1} (A-B))^{-1}  \rangle_r  \leq \frac{1}{1 - \langle A^{-1} (A-B) \rangle_r }\leq \frac{1}{1- \langle A^{-1} \rangle_r \langle A-B \rangle_r} .  $$

\v As we assumed $A$ to be invertible, $B = A(\1_F - A^{-1} (A-B))$ is invertible as well with 

$B^{-1} = (\1_F - A^{-1} (A-B))^{-1} A^{-1}$ and

$$\langle B^{-1} \rangle_r \leq \langle (\1_F - A^{-1} (A-B))^{-1} \rangle_r \langle A^{-1} \rangle_r \leq  \frac{\langle A^{-1} \rangle_r}{1 - \langle A ^{-1} \rangle_r \langle A-B \rangle_r}.$$

\v Finally, $B^{-1} - A^{-1} = B^{-1} (A-B) A^{-1}$ gives

$$\langle B^{-1} - A^{-1} \rangle_r \leq \langle A ^{-1} \rangle_r \langle A-B \rangle_r \langle B^{-1} \rangle_r  \leq \frac{\langle A^{-1} \rangle_r ^2 \langle B -A \rangle_r }{1 - \langle A^{-1} \rangle_r \langle A-B \rangle_r } $$

\v which completes the proof. \hfill \qed

\begin{Theorem}
Let $F$ and $G$ be Fr\'echet spaces, let $U \subset F$ be open. Let $h  \in B^2_r(U,G)$, let $dh(p)$ be bounded-invertible, let $I:  \rightarrow B(G, F)$ be defined by $I(X) = X^{-1}$. Then $J:= I \o dh: V \rightarrow B(G,F)$ is b-differentiable on a possibly smaller open neighborhood $V$ of $p$, and for $q \in V$, the derivative $dJ: V \rightarrow B_r( F , B_r(G,F))$ is given by

$$dJ (v) f = - (dh(v))^{-1} \o dh' (v,f) \o  (dh(v))^{-1} .$$  
\end{Theorem}

\v {\bf Proof.} As the set of bounded-invertible maps $W$ is an open set in $B(F,G)$ according to Theorem \ref{neumann2}, and as $dh: U \rightarrow B(F,G)$ is continuous, $h$ is bounded-invertible in an open neighborhood of $p$. Let $L: B(F,G) \rightarrow B(G,F)$ be the linear map given by $Lf := - A^{-1} f A^{-1} $ which is continuous on $W$ (Theorem \ref{neumann2} again). Now we will show that this expression gives the derivative of the inversion in its domain of definition in $CL(F,G)$. So for $X \in W$,

\bean
I(X) - I(A) - L(X-A) &= X^{-1} - A^{-1} + A ^{-1} (X-A) A^{-1}\\
&= X^{-1} (A-X) A^{-1} + A^{-1} (X-A) A^{-1}\\
&= (-X^{-1} + A^{-1}) (X-A)A^{-1}\\
&= X^{-1} (X-A) A^{-1} (X-A) A^{-1},\\
\eean

\v thus for $X= A+ tB$ we get 

\bean
(I(X) - I(A) - L(X-A))/t &= ((A+tB)^{-1} tB A^{-1} tB A^{-1} )/t\\
&= t (A+tB)^{-1} B A^{-1} B A^{-1} 
\eean

\v As $I$ is continuous according to Theorem \ref{neumann2}, we have $lim_{t \rightarrow 0} (A+tB)^{-1} = A^{-1}$ and

$$lim_{t \rightarrow 0} \frac{I(X) - I(A) - L(X-A)}{t} = lim_{t \rightarrow 0} t (A+tB)^{-1} B A^{-1} B A^{-1}  = 0,$$

\v thus for $dJ$ we get pointwise the expression in the theorem (note that we could not have done the calculation directly in $B_r(F,G)$ as it uses the topological vector space structure). This expression is bounded in $f$ as $h \in B_r^2(F,G)$, and continuous in $v$ by the same reason. This completes the proof. \hfill \qed

\bigskip

\begin{Theorem}[Banach's fix point Theorem]
\label{banach}
Let $(X,d)$ a complete metric space and $f: X \rightarrow X$ a contraction with contraction factor $\rho <1$. Then $f$ has a unique fix point $x_f$ in $X$. It is the limit of the recursive sequence $x_0 \in X$ arbitrary, $x_{n+1} = f(x_n)$. The distance to the solution decreases like

$$d(x_n, x_f) \leq \frac{\rho^n}{1- \rho} d(x_0, x_1) .$$ \hfill \qed
\end{Theorem}

\begin{Theorem}[Left inverse function Theorem between Fr\'echet spaces]
\label{left}
Let $F$ and $G$ be metric Fr\'echet spaces with a metric of the form (\ref{sup}). Let $x_0 \in U$, $U \subset F$ be open, $f: U \rightarrow G$ b-differentiable with radius $R$, let $f'(x_0) $ be injective with $R$-bounded left inverse $L_0 \in B_R(G,F)$. Then there is a real number $r_0 > 0$ such that there is a continuous left inverse of $f \vert_{B(x_0, r_0)}$ (thus in particular it is injective).
\end{Theorem}

\v {\bf Proof.} This is one of the rare occasions where we use local convexity: We choose a convex subset $U_c \subset U $ containing $x_0$. We consider the map $\xi: U \rightarrow B_R(F,F),  \xi: x \mapsto L_0 \o f' (x) - \1_F $ which is bounded-differentiable on $U$ and vanishes at $x_0$. Fix $\rho$ with $0 < \rho < 1$, then by continuity of $df$ we can find a convex $ U' \subset U_c  $ containing $x_0$ such that for $x \in U' $

\bea
\label{vor}
\langle \xi (x) \rangle_R = \langle L_0 f' (x) - \1_F \rangle_R \leq \rho .
\eea

\v thus in $U'$, applying Theorem \ref{convexbound} we get

\bea
\label{abs}
d (L_0 f(x_2) - L_0 f(x_1), x_2 - x_1)  \leq \rho d(x_1, x_2) .
\eea

\v We compute for $x_1, x_2 \in B_r(x_0)$, $ r < R$,

\bean
d_F(x_1, x_2) &= d_F (x_2 - x_1, 0) \\
&\leq d_F(x_2 - x_1, L_0 f (x_2) - L_0 f(x_1)  ) + d_F( L_0 f (x_2) - L_0 f(x_1), 0)\\
&\leq \rho d_F(x_2 - x_1, 0) + \langle L_0 \rangle_R d_G(f(x_2)- f(x_1) , 0), 
\eean

\v therefore

\bea
\label{conse}
\frac{(1 - \rho)}{\langle L_0 \rangle_R } d_F(x_1, x_2 ) \leq d_G(f(x_1), f(x_2)) 
\eea

\v Therefore $f \vert _{B(x_0, r_0)}$ is injective and has an set-theoretic left inverse $\phi := f \vert_{B(x_0, r_0)} ^{-1}$ which takes its image in $U'$ and for which holds $ d_F(\phi (y_1), \phi(y_2)) \leq \frac{\langle L_0 \rangle_R}{1-\rho} d_G(y_1, y_2)  ,$

\v therefore both $f \vert_U$ and its left inverse are bounded, and $\phi$ is continuous. \hfill \qed

\begin{Theorem}[Right inverse function Theorem between Fr\'echet spaces]
Let $F$ and $G$ be metric Fr\'echet spaces with a metric of the form (\ref{sup}). Let $x_0 \in U$, $U \subset F$ be open, $f: U \rightarrow G$ b-smooth of radius $r$ with surjective differential at $x_0$ and bounded right inverse $R_0: U \in B(G,F)$ at $x_0$. Then there is a real number $r_0 > 0$ such that $V:= f(B(x_0, r_0))$ is open in $G$ and there is a b-smooth (of radius $r$) right inverse of $f \vert_{B(x_0, r_0)} \rightarrow V$ (thus in particular it is surjective).
\end{Theorem}
 
\v {\bf Proof.} First, to simplify the notation, by concatenating with translations in $F$ and $G$, let w.r.o.g. $x_0=0, f(x_0) =0$. Rescale the metric $d_G$ such that $\langle f' (0) \rangle_r <1 $. Note that the image $im(R_0)$ of $R_0$ is a Fr\'echet subspace of $F$ as $f' (0)$ is bounded and therefore $R_0$ is bounded from $0$; thus Cauchy sequences in $im(R_0)$ are $R_0$-images of Cauchy sequences in $G$. As in the previous proof, we use local convexity: We choose a convex subset $U_c \subset G $ containing $f(0) =0$. We consider the map $\xi: U_c \rightarrow B(G,G),  \xi: x \mapsto (f' (x) R(0) - \1_G) $ which is bounded-differentiable on $U$ and vanishes at $0$. Fix $\rho$ with $0 < \rho < 1$, then by continuity of $df$ we can find a $ U' \subset U_c \cap B_r (x_0) $ containing $0$ with 

\bea
\label{vor}
\langle \xi \vert_{U'} \rangle_r = \langle f' \vert_{U'} R_0 - \1_G \rangle_r \leq \rho .
\eea

\v Likewise in $U'$, as $f \o R_0$ fixes $0$ and as $(f \o R_0)' = f' \o R_0$, we get

\bea
\label{abs}
d((\1- f \o R_0) y_2, (\1-f\o R_0) y_1) = d (f\o R_0 ( y_2)  - f \o R_0 (y_1) , y_2 - y_1)  \leq \rho d(y_1, y_2) 
\eea

\v for $y_1, y_2 \in U'$, again by applying Theorem \ref{convexbound}. Now we choose a ball radius $ r_0$ such that $\overline{B}(0, r_0) \subset f^{-1} (U') $. Then for $y \in G$ we define $\Phi_y : F \supset \overline{B} (0, r_0) \rightarrow F$ by

$$\Phi_y (x) = x -  (R_0 (f(x) - y)) . $$

\v Note that, up to a possible factor $c>0$ in front of the last bracket, this is exactly the discretization of the real family of operators used in the proof of the Nash-Moser Theorem in \cite{rh}, p. 180. This expression is b-differentiable, and, as $R_0$ is injective, we have $\Phi_y (x) = x \ \Leftrightarrow \ f(x) = y$.

\v The theorem \ref{convexbound} ensures that for $x_1 = R_0 u_1, x_2 = R_0 u_2 \in \overline{B} (0, r_0) \cap im(R_0)  \subset f^{-1} (U')$ we have $d_F(\Phi_y (x_1) , \Phi_y (x_1)) \leq \rho d_F (x_1, x_2)  $ (we use that by surjectivity of $f'(0)$ we know $u_1 = f'(0) x_1$, $u_2 = f' (0) x_2$):

\bean
d(\Phi _y (x_1) , \Phi_y (x_2)) &= d(x_1- R_0 f(x_1), x_2- R_0 f (x_2))\\
&\leq \langle R_0 \rangle_{\infty} d(u_1 - f(R_0 u_1), u_2 - f ( R_0 u_2))\\
&= \langle R_0 \rangle_{\infty} d((\1 - f \o R_0)u_1,( \1 - f \o R_0)u_2 )\\  
&\leq \langle R_0 \rangle_{\infty} \langle \1 - f \o R_0 \rangle_r d(u_1, u_2)\\
&= \langle R_0 \rangle_{\infty} \langle \1 - f \o R_0 \rangle_r d(f'(0)x_1, f' (0)x_2)\\
&\leq \langle R_0 \rangle_{\infty} \langle \1 - f \o R_0 \rangle \langle f' (0) \rangle_r d(x_1, x_2) 
\eean

\v (note that for the third line we need $d(u_1, u_2) < r $ which is guaranteed by $d(u_1, u_2 ) < \langle f' (0) \rangle _{r_0} \cdot r_0$, $r_0 < r$ and $\langle f' (0) \rangle_r <1 $).Thus for appropriate chosen $\rho$, $\Phi_y $ is a contraction for every $y \in G$. We will show now that it also leaves the ball $\overline{B} (0, r_0) \cap im (R_0) $ invariant for appropriately chosen $y$. First of all, by means of linearity of $R_0$ one easily checks that it leaves the subspace $im(R_0)$ invariant. Now let $\b: = \frac{1- \rho}{\langle R_0 \rangle}$ and $r_1:= \beta r_0= \frac{1- \rho}{\langle R_0 \rangle} r_0 $, then for $R_0 (u) =x \in \overline{B} (x_0, r_0) \cap im(R_0)$ and $y \in \overline{B} (f(x_0), r_1)$ we have 

\bean
d_F (0,  \Phi_y (x) ) &= d(R_0 u - R_0 f (R_0 u), R_0 y) \leq \langle R_0 \rangle_{\infty} (d(u- f(R_0 u), 0) + d(y,0))
\eean

\v Therefore if $y \in \overline{B} (f(x_0), r_1)$, then $\Phi_y: \overline{B} (x_0, r_0) \rightarrow \overline{B} (x_0, r_0) $, and we can apply Banach's fix point Theorem \ref{banach} showing that $\Phi_y$ has a unique fixed point in $\overline{B} (x_0, r_0)$, or, equivalently, that for every $y \in \overline{B} (f(x_0), r_1)$ there is a unique $x=: \phi (y) \in \Bb (x_0, r_0)$ with $f(x) = y$. Therefore $\Bb (f(x_0), r_1 ) \subset f(\Bb (x_0, r_0))$.

\bigskip

\v Finally we show that at each $b \in V$ the map $\phi$ is differentiable at $b$ with 

$$\phi '(b) = f' (\phi (b)) ^{-1} $$

\v which implies inductively the smoothness of $\phi'$ on $V$ as a composition of the maps $\phi$, $f'$ and the inversion of bounded-invertible linear maps (cf. Theorem \ref{neumann2}) by help of the chain rule.

\v Differentiability of $\phi$ follows from $\phi (y) := \lim_{k \rightarrow 0} \phi_y^k (x)$ for any $x \in im(R_0)$ and from Theorem \ref{convergence}: To this aim, we observe

$$\frac{d}{dt} \phi_{y + tv}^k (x) = \sum_{i+j = k-1} d \phi_y ^i (\frac{d}{dt} \phi_{y + tv} (\phi _y^j (x)))$$

\v and compute $\frac{d}{dt} \phi_{y + tv} (x) = R_0 v$ (independent of $x$) and $d \phi_y (x) R_0 v = R_0 v - R_0 f' (R_0 v) = R_0 (\1 - f' R_0) (v) $ use again the contraction property of $\phi_y$ and the estimate of $\1- f' \o R_0 $ to prove the local uniform boundedness of this expression. It is easy to see that the family $(\Phi^k_y)'$ is uniformly Lipschitz in $y$. This enables us to use Theorem \ref{convergence}. 

\v Now, knowing differentiability of $\phi$ the equation $\phi '(b) = f' (\phi (b)) ^{-1} $ is an easy consequence from $f \o R = \1_G$ (although it is a nice exercise to show it using the limit above and the Neumann series: As $\Phi_y^k (x) \in im(R_0) $ for all $k$ and all $x \in im(R_0)$, and as $im(R_0)$ is closed, we know that $\phi $ takes its image in $im(R_0)$ as well. Therefore we can consider $\phi_R:= \R_0 ^{-1} \o \Phi: G \supset B(0, r_1) \rightarrow G$. As $R_0 $ is linear and injective, $\phi_R' = R_0 ^{-1} \o \phi'$ which is sufficiently near the identity etc.).     \hfill \qed

\bigskip

\v A corollary of the two preceeding theorems is (as in this case right-inverse and left-inverse coincide)

\begin{Theorem}[Full inverse function Theorem between Fr\'echet spaces]
Let $F$ and $G$ be metric Fr\'echet spaces with the metric being a countable sum of continuous functions of seminorms as given by the equation (\ref{sum}). Let $x_0 \in U$, $U \subset F$ be open, $f: U \rightarrow G$ b-smooth, let $f'(x_0)$ be an isomorphism bounded from $0$. Then there is a real number $r_0 > 0$ such that $V:= f(B(x_0, r_0))$ is open in $G$ and $f \vert_{B(x_0, r_0)} \rightarrow V$ is a diffeomorphism. \hfill \qed
\end{Theorem}

\bigskip

\v As corollaries, we get the inverse function theorems for Fr\'echet manifolds with compatible metrics:

\begin{Theorem}[Left inverse function theorem between metric Fr\'echet manifolds]
Let $M$ resp. $N$ be Fr\'echet manifolds with compatible metrics $d_M$ resp. $d_N$ modelled on metric Fr\'echet spaces $F,G$ with the metric being a countable sum of continuous functions of seminorms as given by the equation (\ref{sum}). Let $x_0 \in U$, $U \subset M$ be open, $f: U \rightarrow N$ bounded-differentiable with $f'': U \times F \rightarrow B(F,G)$ bounded as well in $U$, let $f'(x_0) $ be injective with bounded left inverse $L_0 \in B(G,F)$. Then there is a real number $r_0 > 0$ such that there is a continuous left inverse of $f \vert_{B(x_0, r_0)}$ (thus in particular it is injective). \hfill \qed
\end{Theorem}

\begin{Theorem}[Right inverse function Theorem between metric Fr\'echet manifolds]
Let $M$ resp. $N$ be Fr\'echet manifolds with compatible metrics $d_M$ resp. $d_N$ modelled on metric Fr\'echet spaces $F,G$ with the metric being a countable sum of continuous functions of seminorms as given by the equation (\ref{sum}). Let $x_0 \in U$, $U \subset M$ be open, $f: U \rightarrow G$ bounded-differentiable with surjective differential at $x_0$ and bounded right inverse $R_0: U \in B(G,F)$ at $x_0$ and with $f'': U \times F \rightarrow B(F,G)$ bounded as well in $U$. Then there is a real number $r_0 > 0$ such that $V:= f(B(x_0, r_0))$ is open in $G$ and there is a smooth right inverse of $f \vert_{B(x_0, r_0)} \rightarrow V$ (thus in particular it is surjective). \hfill \qed
\end{Theorem}

\begin{Theorem}[Full inverse function theorem between metric Fr\'echet manifolds]
Let $M$ resp. $N$ be Fr\'echet manifolds with compatible metrics $d_M$ resp. $d_N$ modelled on metric Fr\'echet spaces $F,G$ with a metric of the form (\ref{sup}). Let $U \subset M$ be an open set contained in a chart, $f: U \rightarrow N$ bounded-differentiable and $f'': U \times F \rightarrow B(F,G) $ bounded as well, everything w.r.t. $d_M$ and $d_N$. Let $f(x_0) $ be an isomorphism bounded from $0$. Then there is a real number $r_0 > 0$ such that $V:= f(B(x_0, r_0))$ is open in $N$ and $f \vert_{B(x_0, r_0)} \rightarrow V$ is a diffeomorphism. \hfill \qed  
\end{Theorem}

\section{Fr\'echet manifolds as generalized length spaces}

\v Now let us equip the Fr\'echet spaces and manifolds with structures carrying still a little bit more information trying to give them a notion of curve length. We recall the notion of {\bf length spaces}, adding two generalizations:

\begin{Definition}
A {\bf length structure} on a topological space $X$ is a subset $C$ of $\bigcup_{(a,b \in \R^2)} C^0 ([a,b], X) $ closed under concatenation of paths, restrictions to subsets and linear reparametrizations, and a function $L: C \rightarrow \R \cup \{ \infty \} $ which is additive under concatenation of curves, invariant under linear homeomorphisms and with the additional properties that for a curve $c: [a,b] \rightarrow X$ with $L(c) < \infty$ the function $l: [a,b] \rightarrow \R^+ , t \mapsto L (c \vert_{[a,r]})$ is continuous, and that for every point $x \in X$ and every open neighborhood $U$ of $x$ there is an $\e >0$ such the length of any curve beginning at $x$ and ending in the complement of $U$ is greater than $\e$. A {\bf generalized length structure of first kind} is a subset and a function satisfying all of the above, except that the function does not have to be additive under concatenation of curves. A {\bf generalized length structure of second kind} is a subset and a function satisfying all of the above, except that the function does not have to be invariant under reparametrizations of curves. A {\bf length space} is a metric space such that the distance between two points $p,q$ coincides with the infimum of length of curves joining $p$ and $q$.
\end{Definition}

\v If $X$ is a metric space, there is always the (in the category of metric spaces and isometries) natural choice of a length structure $C_0:= \bigcup_{a,b \in \R^2} C^0 ([a,b], X), L_0 := sup_{P \in \mathcal{P}} \sum_{i=1}^n d(c(t_{i+1}), c(t_i))$. If a curve $k: [a,b] \rightarrow X$ is Lipschitz, i.e. if $lip_c := sup_{s,t \in [a,b], s \neq t} d(k(s), k(t))/\vert s-t \vert < \infty$, then $L_0 (k)= \int_a^b lip_t dt$, where $lip_t := \lim_{\e \rightarrow 0} lip (c \vert_{[t- \e, t + \e]})$ (this limit exists as its argument is monotonely decreasing in $\e$ and bounded by $0$). The definition of $L_0$, due to Gromov, is an extension of the notion of curve length for piecewise differentiable curves in Banach manifolds. However, the following corollary from \ref{fine} shows that this is not the appropriate notion of length as we would have too few rectifiable curves, i.e. curves of finite length, as these would not even include all straight lines:

\begin{Corollary} 
Let $g: t \mapsto tv$ be a straight line in $\G (\pi)$ with a metric of the form \ref{sum} or \ref{sup}, then $L_0 (g) < \infty \Leftrightarrow \exists M \in \R : \vert \vert v \vert \vert_i \leq M \forall i $. 
\end{Corollary}

\v So we try to construct at least {\em generalized} length structures for interesting spaces of curves. To this purpose we first prove a theorem about arc-length parametrization of differentiable curves:

\begin{Theorem}
Given a $C^1$ curve $c: [a,b] \rightarrow F$ in a Fr\'echet space with $\dot{c} (t) \neq 0$, we can always find a reparametrization $\tilde{c}: [0, T] \rightarrow X$ with $d_F(\dot{c} (t), 0) =1$.  
\end{Theorem}

\v {\bf Proof.} This follows from continuity of the first Minkowski functional and the multiplicative inversion in $\R^+$. \hfill \qed

\bigskip

\v We define the {\bf metric length $l$} of a continuous curve $c: \In \rightarrow F$ by 

$$l(c) := \int_{\In} \sum \a_n \Phi (sup_{r \in \R \setminus \{ 0 \} } \vert \vert \dot{c}(t) (r \partial_t)  \vert \vert_n^M )  dt$$ 

\v where the $\vert \vert \cdot \vert \vert^M_k$ are the Minkowski functionals. 

\v We define the {\bf smooth length $L$} of a smooth curve by the canonical embedding $H: C^{\infty } (\In, \G (\pi)\rightarrow \G (\In \times \pi)  $ and $L(H) := \langle H \rangle $ as in the section of Fr\'echet manifolds. 

\v Obviously, the smooth length gives rise to the Fr\'echet metric already used, but its definition needs more than the metric alone, but the whole series of seminorms. The metric length, in turn, is defined for a wider class of curves, but it induces a new metric. It seems most approriate to define lengths of smooth curves {\em first}, by $L$, and then define the metric as infimum of curve length. We want to compare $l$ and $L$ which differ basically by the order of $\Phi$ and the integral along the curve:

\begin{Theorem}
Let $M$ be a metric Fr\'echet manifold. $(C^{\infty}( \In, M), L)$ is a generalized length structure of first kind, satisfying only an inequality $L(c_1 \o c_2) \leq L(c_1) + L(c_2)$. $(C^0 (\In, M), l)$ is a generalized length structure of second kind. We have the inequality $l(c) \leq L(c)$. 
\end{Theorem}

\v {\bf Proof.} We only show $\langle I \rangle \geq l(I)$. It is enough to show the inequality in

$$\sum 2^{-n} \Phi (a_n \int_0^1 f(t) dt  ) \geq \int_0^1 \sum 2^{-n } \Phi (a_nf(t) ) dt =  \sum 2^{-n } \int_0^1 \Phi (a_nf(t) ) dt    ,$$

\v as the right hand side is $\geq \int [...] dt$ by the inequality after the definition of length and scalar-boundedness. So if we absorb $a_m$ into $f$, it is enough to show that for an arbitrary positive function $S$, we have 

\bea
\label{hotte}
\Phi (\int_0^1 f(t) dt ) \geq \int_0^1 \Phi (f(t)) dt .
\eea

\v But this is clear from the concavity of $\Phi$ and the fact that step functions are dense in $L^1(\In)$: Let $T$ be a step function on the unit interval with values $T_1, ... T_n$. The integral of $T$ is a convex combination of $T_1,... T_n$, and a concave function maps convex sets into convex sets, explicitly: $\Phi (tx + (1-t) y) \geq t\Phi(x) + (1-t) \Phi (y) $, thus by induction $\Phi (\sum t_i x_i ) \geq \sum t_i \Phi (x_i)$ if $\sum t_i =1$. Therefore Equation (\ref{hotte}) holds, and we are done. \hfill \qed

\bigskip

\v Some important questions remain open: Is there a metric such that Gromov's curve length is finite for all smooth curves? Is there a proper length structure (no generalized one) only defined by means of the metric for which all smooth curves are rectifiable? And if yes, is there one inducing the metric it started with? This shall be object of further research.

\end{document}